\documentclass[12pt]{article}

\usepackage{amssymb,amsthm,amsmath,eucal,latexsym}
\usepackage{color,caption,comment,subcaption}
\usepackage{multirow}
\usepackage{listings}
\usepackage{setspace}

\usepackage{changepage}
\usepackage[linesnumbered,lined,boxed,commentsnumbered]{algorithm2e}
\usepackage{verbatim}
\usepackage{tikz}
\usepackage{xcolor}
\usepackage{algpseudocode}

\newtheorem{theorem}{Theorem}[section]
\newtheorem{corollary}[theorem]{Corollary}
\newtheorem{lemma}[theorem]{Lemma}
\newtheorem{conjecture}[theorem]{Conjecture}
\newtheorem{definition}[theorem]{Definition}

\usepackage[utf8]{inputenc}

\usepackage{enumitem, graphicx, hyperref}
\allowdisplaybreaks

\title{Extending Graph Burning to Hypergraphs}
\author{Andrea C.~Burgess\thanks{Department of Mathematics and Statistics, University of New Brunswick, Saint John, NB, E2L~4L5, Canada. \texttt{andrea.burgess@unb.ca}}, Caleb W.~Jones\thanks{Department of Mathematics, Toronto Metropolitan University, Toronto, ON, M5B~2K3, Canada. \texttt{caleb.w.jones@torontomu.ca}}, David A.~Pike\thanks{Department of Mathematics and Statistics, Memorial University of Newfoundland, St.~John's, NL, A1C~5S7, Canada. \texttt{dapike@mun.ca}}}

\begin{document}

\date{}
\maketitle


\begin{abstract}
Graph burning is a round-based game or process that discretely models the spread of influence throughout a network. We introduce a generalization of graph burning which applies to hypergraphs, as well as a variant called ``lazy'' hypergraph burning. Interestingly, lazily burning a graph is trivial, while lazily burning a hypergraph can be quite complicated. Moreover, the lazy burning model is a useful tool for analyzing the round-based model. One of our key results is that arbitrary hypergraphs do not satisfy a bound analogous to the one in the Burning Number Conjecture for graphs. We also obtain bounds on the burning number and lazy burning number of a hypergraph in terms of its parameters, and present several open problems in the field of (lazy) hypergraph burning.
\end{abstract}


\section{Introduction}
\label{intro}

Graph burning is a single-player game played on finite, simple, undirected graphs over a discrete sequence of rounds. The player, whom we call the \textit{arsonist}, attempts to set fire to every vertex of the graph in as short a time as possible. The arsonist manually sets fire to vertices, and the fire also spreads or \textit{propagates} from burned vertices to adjacent unburned vertices in each round. Once a vertex is set on fire, it remains on fire until the end of the game. Of course, the game ends when every vertex in the graph is on fire.

A problem equivalent to a special case of graph burning was posed in 1992 by Brandenburg and Scott as an internal problem at Intel \cite{VERY_first_ref}. It modelled the transfer of information through processors, and was equivalent to burning an $n$-cube. Independently in 2014, Bonato, Janssen, and Roshanbin \cite{first_paper} introduced graph burning as a combinatorial process which could be applied to any graph. They introduced the term ``graph burning'' for the first time, and developed much of the theory. 

We now describe the rules for graph burning in detail. Rounds are indexed by $\mathbb{N}$ starting at one. Denote the set of vertices that are on fire at the end of round $r$ by $F_r$, and set $F_0=\emptyset$. During each round $r\geq 1$, the following two things happen simultaneously. 

\begin{itemize}
\item For each $v\in V(G)\setminus F_{r-1}$, if there is $u\in F_{r-1}$ such that $vu\in E(G)$ then $v$ catches fire. 
\item The arsonist chooses a vertex $u_r \in V(G)\setminus F_{r-1}$ and sets it on fire ($u_r$ is called a \textit{source}).
\end{itemize}

In round 1, no vertices catch fire due to propagation, and the arsonist chooses the first source. Note that the arsonist may choose as a source a vertex that also catches fire due to propagation that round, although it is never advantageous to do so. We will call such a source \textit{redundant}. In the last round, the arsonist may have no choice but to choose a redundant source. A sequence of sources $(u_1,u_2,\ldots,u_k)$ is called \emph{valid} if $u_r\notin F_{r-1}$ for each $r$, and otherwise it is called \emph{non-valid}. Observe that a non-valid sequence of sources prompts the arsonist at least once to burn as a source a vertex that was on fire at the end of the previous round, which is not allowed.

A valid sequence of sources $(u_1,u_2,\ldots,u_k)$ that leaves the graph completely burned when the arsonist burns $u_i$ in round $i$ is called a \textit{burning sequence}. A burning sequence of minimum possible length is called \textit{optimal}. A graph may have many different optimal burning sequences. Given a graph $G$, the \textit{burning number} of $G$, denoted $b(G)$, is a measure of how fast the fire can possibly spread to all vertices of $G$. In particular, it is the earliest round at which $G$ could possibly be completely burned. Note that $b(G)$ is also equal to the length of an optimal burning sequence since the number of rounds coincides with the number of sources that are chosen. 

The following conjecture is perhaps the deepest open problem in the field of graph burning. The two theorems that follow are fundamental results on the burning number. 

\begin{conjecture}
\label{burn_num_conj}
\emph{\textbf{(}Burning Number Conjecture, \cite{first_paper}\textbf{)}} For a connected graph $G$ of order $n$, $b(G) \leq \lceil \sqrt{n}\ \rceil$. 
\end{conjecture}

\begin{theorem}
\label{tree_red_thm}
\emph{\textbf{(}Tree Reduction Theorem, \cite{first_paper}\textbf{)}} For a graph $G$, $$b(G)=\min \{b(T) \mid T \text{ is a spanning subtree of } G \}.$$
\end{theorem}

\begin{theorem}
\label{path_thm}
\emph{\textbf{(}\cite{first_paper}\textbf{)}} For a path $P_n$ on $n$ vertices, $b(P_n)=\lceil \sqrt{n}\ \rceil$.
\end{theorem}

Conjecture \ref{burn_num_conj} and Theorem \ref{path_thm} together would imply that paths are in a sense the ``hardest'' graphs to burn. Theorem \ref{tree_red_thm} is a useful tool for making progress towards the Burning Number Conjecture -- if we are trying to find an upper bound on $b(G)$ in terms of $|V(G)|$ for an arbitrary graph $G$, we may assume $G$ is a tree.

Recent progress towards the Burning Number Conjecture has yielded the following two results, the second proving that the Burning Number Conjecture holds asymptotically.

\begin{theorem}
\emph{\textbf{(}\cite{new_paper}\textbf{)}} For a connected graph $G$ of order $n$, $b(G) \leq 1+\left\lceil \sqrt{\frac{4n}{3}}\ \right\rceil$. 
\end{theorem}

\begin{theorem}
\emph{\textbf{(}\cite{asymptotically_paper}\textbf{)}} For a connected graph $G$ of order $n$, $b(G) \leq (1+o(1))\sqrt{n}$.
\end{theorem}

See \cite{bounds,Bonato_summary,how_to} for more on graph burning.

We now give a brief review of some definitions from hypergraph theory. For more information on hypergraphs, see \cite{bahm_sanj,hypergraph_text,conn_in_hyp}. A \emph{hypergraph} $H$ is an ordered pair $(V,E)$ where $V$ and $E$ are disjoint finite sets, $V\neq\emptyset$, and each element of $E$ is a subset of $V$. The elements of $V=V(H)$ are called \emph{vertices}, and the elements of $E=E(H)$ are called \emph{edges} or \emph{hyperedges}. Informally, a hypergraph is like a graph, except edges can now contain any number of vertices, not just two (notice that every graph is also a hypergraph). If two vertices belong to a common edge, they are called \emph{adjacent}. The \emph{degree} of a vertex $v$ is the number of edges containing $v$. An alternating sequence of vertices and edges $v_1,e_1,v_2,e_2,\ldots,v_n,e_n,v_{n+1}$ is called a \emph{path} if $v_1,v_2,\ldots,v_{n+1}$ are all distinct vertices, $e_1,e_2,\ldots,e_n$ are all distinct edges, and $v_i,v_{i+1}\in e_i$ for each $i\in\{1,2,\ldots,n\}$. The \emph{length} of the path is $n$, and the path \emph{connects} $v_1$ and $v_{n+1}$. A hypergraph is \emph{connected} if, for any two vertices $x$ and $y$, there is a path that connects $x$ and $y$. Otherwise, the hypergraph is \emph{disconnected}. If a vertex does not belong to any edge, it is called \emph{isolated}. A hypergraph is \emph{linear} if any two distinct edges intersect in at most one vertex. A hypergraph is $k$\emph{-uniform} if every edge contains exactly $k$ vertices. Two edges $e_1$ and $e_2$ are \emph{parallel} if they contain exactly the same vertices. The number of edges parallel to some edge $e$, including $e$ itself, is the \emph{multiplicity} of $e$. A hypergraph is called \emph{simple} if no edge has multiplicity greater than one, and no edge contains one or fewer vertices. The $2$-\emph{section} of a hypergraph $H$ is the graph $G$ on the same vertex set whose edges are precisely those of the form $\{u,v\}$ such that $\{u,v\}\subseteq e$ for some $e\in E(H)$.

In this paper we introduce a generalization of graph burning that is played on hypergraphs. The new game should look much the same -- each round, the arsonist burns a vertex, and the fire spreads to other ``nearby'' vertices based on some propagation rule. This rule should ensure that the game on hypergraphs reduces to the original game when each edge of the hypergraph contains exactly two vertices (i.e.~when we play it on a graph). The most obvious rule is that fire spreads from burned vertices to adjacent unburned vertices; however this is equivalent to burning the $2$-section of the hypergraph, so this is essentially a special case of graph burning.



We therefore formulate the following rule for how the fire propagates. Fire spreads to a vertex $v$ in round $r$ if and only if there is a non-singleton edge $\{v,u_1,\ldots,u_k\}$ such that each of $u_1,u_2,\ldots,u_k$ was on fire at the end of round $r-1$. Clearly this reduces to the original game when the hypergraph has only edges of size two. Rounds are indexed by $\mathbb{N}$ starting at one, they have the same structure as in graph burning, and the definitions of a \textit{source}, \textit{redundant source}, \textit{burning sequence}, and \textit{burning number} are all analogous to those in graph burning.

We also introduce an alternate set of rules for how the arsonist burns a hypergraph. Suppose the arsonist is very lazy, and does not wish to be present while the hypergraph burns. In particular, they wish to set fire to a select few vertices simultaneously in such a way that the hypergraph is eventually completely burned through subsequent propagation. Of course, the arsonist wishes to set fire to as few vertices as possible while still ensuring that the hypergraph becomes completely burned through propagation. We call the set of vertices the arsonist initially sets fire to a \textit{lazy burning set}. The size of a smallest or \textit{optimal} lazy burning set for a hypergraph $H$ is called the \textit{lazy burning number} of $H$, denoted $b_L(H)$. 

The lazy burning game on graphs is trivial. The arsonist must simply set fire to exactly one vertex in each connected component of the graph to achieve a minimum lazy burning set. However, lazily burning a hypergraph is much more interesting. Furthermore, we will show that the lazy burning model is a useful tool for analyzing the round-based model of hypergraph burning. 

We observe that lazy hypergraph burning already exists in the literature, having been introduced in \cite{bootstrap_paper_1} under the name $\mathcal{H}$\emph{-bootstrap percolation}. The existing results on $\mathcal{H}$-bootstrap percolation are mostly probabilistic or extremal in nature, or apply to specific families of hypergraphs such as hypercubes. Our results on lazy hypergraph burning take a different approach than those in the literature, as they are deterministic in nature, and focus on the connection between the lazy and round-based versions of the game. See also \cite{another_crit_prob,threshhold_paper,threshhold_paper_1,crit_prob,crit_prob_again} for a review of $\mathcal{H}$-bootstrap percolation.

For the remainder of this paper we mainly focus on bounding the burning number and lazy burning number of a hypergraph in terms of its parameters. This is achieved (with varying degrees of success) for arbitrary hypergraphs, tight $3$-uniform paths, disconnected hypergraphs, and subhypergraphs. We also have a number of supplementary results. For example, we prove that there is no analogue to the Burning Number Conjecture for either round-based or lazy hypergraph burning.

\section{General Results and Bounds}
\label{general_results}

The burning and lazy burning numbers of a hypergraph can be bounded above and below by simple hypergraph parameters. The following results establish these bounds individually, and they are combined in Theorem \ref{cool_eqtn}.

Given a hypergraph $H$, let $\mathcal{E}(H)$ be the number of edges in $H$ that are not singleton, empty, or duplicate edges. Note that if an edge has multiplicity greater than one, then we choose one instance of the edge to be the ``original,'' and the rest are ``duplicates.'' Thus, an edge with multiplicity greater than one contributes $1$ to the sum $\mathcal{E}(H)$ (provided that it is not a singleton or empty edge). 

\begin{theorem}
Let $H$ be a hypergraph, and let $\mathcal{E}(H)$ be the number of edges in $H$ that are not singleton, empty, or duplicate edges. Then $|V(H)|-\mathcal{E}(H) \leq b_L(H)$.
\end{theorem}

\begin{proof}
Let $S$ be a lazy burning set for $H$. Denote by $z_S$ the number of vertices that become burned throughout the game through propagation, so $z_S+|S|=|V(H)|$. Observe that each edge (that is not a singleton, empty, or duplicate edge) may ``cause'' fire to spread to at most one vertex throughout the lazy burning game. Thus, $z_S\leq \mathcal{E}(H)$. But then $|S|=|V(H)|-z_S\geq|V(H)|-\mathcal{E}(H)$. In particular, if $S$ is a minimum lazy burning set we get $b_L(H)=|S|\geq|V(H)|-\mathcal{E}(H)$.
\end{proof}

Of course, if $H$ is simple then $\mathcal{E}(H)=|E(H)|$. We therefore have the following corollary.

\begin{corollary}
\label{bl_lower_bound}
For any simple hypergraph $H$, $|V(H)|-|E(H)|\leq b_L(H)$. 
\end{corollary}

\begin{theorem}
\label{loose_ineq}
$b_L(H) \leq b(H)$ for all hypergraphs $H$.
\end{theorem}

\begin{proof}
Given any burning sequence $(u_1,u_2,...,u_k)$ for $H$, $\{u_1,u_2,...,u_k\}$ is a lazy burning set for $H$, and hence $b_L(H)\leq k$. Choose an optimal burning sequence $(u_1,u_2,...,u_{b(H)})$ for $H$. Then $b_L(H)\leq b(H)$.
\end{proof}

\begin{corollary}
Let $H$ be a hypergraph, and let $\mathcal{E}(H)$ be the number of edges in $H$ that are not singleton, empty, or duplicate edges. Then $|V(H)|-\mathcal{E}(H)\leq b(H)$.
\end{corollary}

\begin{corollary}
\label{edge_idea}
For any simple hypergraph $H$, $|V(H)|-|E(H)|\leq b(H)$. 
\end{corollary}



\begin{figure}[h]
\centering
\begin{tikzpicture}[scale=1.1]

\node (x) at (0,1) {};
\fill [fill=black] (x) circle (0.105) node [below] {$x$};
\node (y) at (1,0) {};
\fill [fill=black] (y) circle (0.105) node [above] {$y$};
\node (z) at (0,-1) {};
\fill [fill=black] (z) circle (0.105) node [above] {$z$};
\node (w) at (-1,0) {};
\fill [fill=black] (w) circle (0.105) node [right] {$w$};

\draw [line width=0.25mm] [black] [rounded corners=0.7cm] (0.25,1.5)--(0.25,-1.5)--(-1.5,0)--cycle;
\end{tikzpicture}

\caption{An example of a hypergraph $H$ where all the bounds in Corollary \ref{bl_lower_bound}, Theorem \ref{loose_ineq}, and Corollary \ref{edge_idea} are tight simultaneously. Note that $|E(H)|=1$ since singleton edges are not allowed in simple hypergraphs.}
\label{disconn1}
\end{figure}
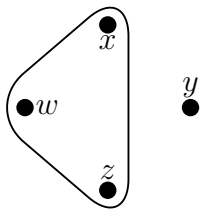

The bounds in Corollary \ref{bl_lower_bound}, Theorem \ref{loose_ineq}, and Corollary \ref{edge_idea} are tight. Consider the simple hypergraph $H$ in Figure \ref{disconn1}. A minimum lazy burning set is $\{x,y,z\}$, an optimal burning sequence is $(x,z,y)$, and $|V(H)|-|E(H)|=3$. Hence, $|V(H)|-|E(H)|=b_L(H)=b(H)$, so the aforementioned bounds are all tight simultaneously. This example generalizes to an infinite family of hypergraphs that show tightness; one may construct such a hypergraph $H$ with $V(H)=\{v_1,\ldots,v_n\}$ and $E(H)=\big{\{}\{v_1,\ldots,v_{n-1}\}\big{\}}$, where $n\geq 3$. Then a minimum lazy burning set is $\{v_2,\ldots,v_n\}$, an optimal burning sequence is $(v_2,\ldots,v_n)$, and $|V(H)|-|E(H)|=n-1$.


The bound in Theorem \ref{loose_ineq} can be improved to a strict inequality if $H$ has no isolated vertices (or if a less strict but harder to discern condition is met; see Theorem \ref{strict_ineq}).

\begin{theorem}
\label{strict_ineq}
If there is an optimal burning sequence $(u_1,u_2,\ldots,u_{b(H)-1}, u_{b(H)})$ in $H$ such that the last source $u_{b(H)}$ is not an isolated vertex, then $b_L(H)<b(H)$.
\end{theorem}

\begin{proof}
We will show that $\{u_1,u_2,\ldots,u_{b(H)-1}\}$ is a lazy burning set for $H$. Since $u_{b(H)}$ is the final source, burning $u_1,u_2,\ldots,u_{b(H)-1}$ one-by-one (as in the original game) or simultaneously (as a lazy burning set) will eventually result in all of $V(H)\setminus \{u_{b(H)}\}$ being burned. So, let us burn each vertex in $\{u_1,u_2,\ldots,u_{b(H)-1}\}$ simultaneously. All we need to show is that $u_{b(H)}$ will eventually burn through propagation. If $u_{b(H)}$ is a redundant source then clearly fire will propagate to $u_{b(H)}$. Otherwise, $u_{b(H)}$ is not a redundant source. But eventually all of $V(H)\setminus \{u_{b(H)}\}$ will burn through propagation, and $u_{b(H)}$ belongs to an edge since it is not isolated. All other vertices in the edge containing $u_{b(H)}$ are on fire, and thus $u_{b(H)}$ will catch on fire. 
\end{proof}

\begin{corollary}
\label{cool_cor_1}
If $|V(H)|\geq 2$ and $H$ has no isolated vertices then $b_L(H)<b(H)$.
\end{corollary}

\begin{corollary}
\label{cool_cor_2}
If $|V(H)|\geq 2$ and $H$ is connected then $b_L(H)<b(H)$.
\end{corollary}

If $V(H)=\{v_1,\ldots,v_n\}$ and $E(H)=\big{\{} \{v_1,\dots,v_k\},$ $\{v_k,\ldots,v_{2k-1}\},\ldots,\{v_{n-k+1},\ldots,v_n\} \big{\}}$, then $H$ is called a $k$\emph{-uniform loose path}. Informally, a \emph{loose path} is any hypergraph that can be created from a $k$-uniform loose path by deleting vertices of degree one while ensuring no edge becomes a singleton. 

\begin{figure}[h]
\centering
\begin{tikzpicture}[scale=1.1]

\node (1) at (0,1) {};
\fill [fill=black] (1) circle (0.08) node [below] {};
\node (2) at (0.71,0.71) {};
\fill [fill=black] (2) circle (0.08) node [below left] {};
\node (3) at (1,0) {};
\fill [fill=black] (3) circle (0.08) node [left] {};
\node (4) at (0.71,-0.71) {};
\fill [fill=black] (4) circle (0.08) node [above left] {};
\node (5) at (0,-1) {};
\fill [fill=black] (5) circle (0.08) node [above] {};
\node (6) at (-0.71,-0.71) {};
\fill [fill=black] (6) circle (0.08) node [above right] {};
\node (7) at (-1,0) {};
\fill [fill=black] (7) circle (0.08) node [right] {};
\node (8) at (-0.71,0.71) {};
\fill [fill=black] (8) circle (0.08) node [below right] {};

\draw [line width=0.25mm] [black] (0,0) ellipse (1.25 cm and 1.25cm);
\end{tikzpicture}

\caption{An example where all the bounds in Theorem~\ref{cool_eqtn} are tight simultaneously.}
\label{disconn2}
\end{figure}
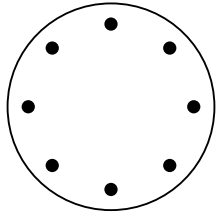

A hypergraph $H$ for which $b_L(H)=b(H)-1$ can be seen in Figure \ref{disconn2}, so the bounds from Theorem \ref{strict_ineq} and its two corollaries are tight. There is indeed an infinite family of hypergraphs that exhibit the tightness of these bounds -- the family of loose paths with minimum edge size three. Denote the edges of a loose path by $e_1, e_2, \ldots, e_m$ such that $e_1$ and $e_m$ are the edges that contain exactly one vertex of degree two, and such that $e_i$ and $e_{i+1}$ share a vertex for each $i\in\{1,2,\ldots,m-1\}$. One may construct an optimal burning sequence in a loose path (with minimum edge size three) by burning all the degree-one vertices in $e_1$ as sources, followed by all the degree-one vertices in $e_2$, and so on. Indeed, by following this process, every degree-one vertex will be a source in the burning sequence, and the last source will be redundant. Furthermore, one may construct a minimum lazy burning set by taking the vertices in an optimal burning sequence as a set and deleting any one vertex. Hence, each loose path $H$ with minimum edge size three has $b_L(H)=b(H)-1$.

A set of vertices $\{x_1,x_2,\ldots,x_k\}\subseteq V(H)$ is called \emph{independent} if there is no edge $e$ in $H$ such that $e\subseteq\{x_1,x_2,\ldots,x_k\}$. The size of a largest independent set in $H$ is the \emph{independence number} of $H$, denoted $\alpha(H)$. We now consider upper bounds on $b_L(H)$ and $b(H)$ which make use of $\alpha(H)$. 

\begin{lemma}
\label{opt_ind_set_lazy}
Any optimal lazy burning set in a hypergraph $H$ is an independent set.
\end{lemma}

\begin{proof}
Suppose $S$ is an optimal lazy burning set in $H$ that is not independent. Then there is some edge $e=\{x_1,x_2,\ldots,x_k\} \subseteq S$ in $H$. But then $S \setminus \{x_k\}$ is a smaller lazy burning set in $H$, which contradicts $S$ being optimal. 
\end{proof}

\begin{lemma}
\label{lazy_lem_2}
For a hypergraph $H$, any independent set of vertices of size $\alpha(H)$ is a lazy burning set for $H$.
\end{lemma}

\begin{proof}
Let $S\subseteq V(H)$ be an independent set with $|S|=\alpha(H)$. Consider any vertex $v \in V(H) \setminus S$. Clearly $\{v\}\cup S$ is not an independent set by the maximality of $S$. Thus, there is an edge $\{v,x_1,x_2,\ldots,x_k\}\subseteq \{v\}\cup S$ in $H$. But then if we set $S$ on fire, each of $x_1,x_2,\ldots,x_k$ will be set on fire, and thus the fire will spread to $v$ through propagation. The vertex $v$ was arbitrarily chosen, so each vertex in $V(H) \setminus S$ will burn through propagation. Therefore $S$ is a lazy burning set for $H$.
\end{proof}

Note that both Lemma \ref{opt_ind_set_lazy} and Lemma \ref{lazy_lem_2} imply Theorem \ref{lazy_lem_1}.

\begin{theorem}
\label{lazy_lem_1}
$b_L(H) \leq \alpha(H)$ for all hypergraphs $H$.
\end{theorem}

The bound in Theorem \ref{lazy_lem_1} is tight; see Figure \ref{disconn1} for an example. The hypergraph pictured has independence number and lazy burning number three, as $\{x,y,z\}$ is both a maximum independent set and a minimum lazy burning set. Indeed, there is an infinite family of hypergraphs that exhibit the tightness of Theorem \ref{lazy_lem_1}. One may construct such a hypergraph $H$ with $V(H)=\{v_1,\ldots,v_n\}$ and $E(H)=\big{\{}\{v_1,\ldots,v_{n-1}\}\big{\}}$, where $n\geq 3$. Then $\{v_2,\ldots,v_n\}$ is both a maximum independent set and an minimum lazy burning set for $H$, so $b_L(H) = \alpha(H)$.


\begin{theorem}
\label{b_alphaplus1}
$b(H)\leq \alpha(H)+1$ for all hypergraphs $H$. 
\end{theorem}

\begin{proof}
Let $S=\{u_1,u_2,\ldots,u_{\alpha(H)}\}$ be a maximum independent set in $H$. Let the arsonist burn the $u_i$ in order as a burning sequence. If at some round the arsonist expects to burn some $u_j$ as a source, but it is already on fire, then the arsonist skips $u_j$ (which only shortens the burning sequence). Thus, at the end of some round $r\leq \alpha(H)$, each vertex in $S$ is on fire (and possibly some others). Assume that at the end of round $r$, $H$ is not fully burned (otherwise, we are done). We claim that in the following round $r+1$ the rest of $H$ will burn through propagation (and the arsonist chooses a redundant source).

Consider any vertex $v\in V(H)$ that was not on fire at the end of round $r$. Suppose $v$ does not catch fire in round $r+1$. Then, at the end of round $r$, there was no edge $e$ containing $v$ such that all of $e\setminus\{ v\}$ was on fire. In particular, since $S$ is completely burned at the end of round $r$, there is no edge $e$ containing $v$ such that $e\setminus\{ v\} \subseteq S$. But then $S\cup\{ v\}$ is an independent set in $H$ that is strictly larger than $S$, which is a contradiction. Therefore, each vertex of $H$ is on fire at the end of round $r+1$, so $b(H)\leq r+1 \leq \alpha(H)+1$.
\end{proof}


The bound in Theorem \ref{b_alphaplus1} is also tight; see Figure \ref{disconn2} for an example. Again, there is an infinite family of hypergraphs that exhibit the tightness of this bound. Consider the family of hypergraphs that consist only of disjoint non-singleton edges. In such a hypergraph, one may construct an optimal burning sequence by burning the vertices of a maximum independent set in any order, with an additional redundant source in the final round. Hence, each hypergraph $H$ in this family has $b(H)= \alpha(H)+1$.

The difference between $b_L(H)$ and $\alpha(H)$ (respectively $b(H)$ and $\alpha(H)+1$) can also be arbitrarily large. Consider the family of hypergraphs $H$ with $V(H)=\{v_1,\ldots,v_n\}$ and $E(H)=\big{\{}\{v_1,v_2,v_3\}, \{v_1,v_2,v_4\},\ldots, \{v_1,v_2,v_n\}\big{\}}$. Such a hypergraph has $b_L(H)=2$, $b(H)=3$, and $\alpha(H)=n-1$.






The following result is immediate due to Corollary \ref{cool_cor_1} and Theorem \ref{b_alphaplus1}.

\begin{corollary}
\label{alpha_implies_alphaplusone}
Let $H$ be a hypergraph with no isolated vertices. If $b_L(H)=\alpha(H)$ then $b(H)=\alpha(H)+1$.
\end{corollary}

We are thus far unaware of any conditions on $H$ that are sufficient for concluding $b_L(H)=\alpha(H)$. In particular, we ask whether or not the converse of Corollary \ref{alpha_implies_alphaplusone} is true, as this would be one such sufficient condition. 

Finally, by combining the bounds in Corollary \ref{bl_lower_bound}, Corollary \ref{cool_cor_1}, and Theorem \ref{b_alphaplus1}, we get the series of inequalities in the following result.  

\begin{theorem}
\label{cool_eqtn}
Let $H$ be a simple hypergraph with no isolated vertices. Then $$|V(H)|-|E(H)|\leq b_L(H)<b(H)\leq \alpha(H)+1.$$
\end{theorem}

Each inequality in Theorem \ref{cool_eqtn} is tight, and Figure \ref{disconn2} shows an example of a hypergraph $H$ where all of them are tight simultaneously. It has $7=|V(H)|-|E(H)|=b_L(H)$ and $8=b(H)=\alpha(H)+1$. Of course, this example can be expanded to an infinite family of hypergraphs that exhibit tightness for each inequality in Theorem \ref{cool_eqtn} simultaneously. Simply consider the family of hypergraphs that consist of a single edge containing all of their vertices (excluding the hypergraph which is a single vertex in an edge).


Considering the inequalities in Theorem~\ref{cool_eqtn} further, we have already shown that the difference between $\alpha(H)+1$ and $b(H)$ can be arbitrarily large.  Additionally it is easy to see that $b_L(H)-(|V(H)|-|E(H)|)$ can also be arbitrarily large, by noting that $b_L(H)>0$ for every hypergraph $H$, but, for instance, taking $H$ to be the $k$-uniform complete hypergraph of order $n$, we have that $|V(H)|-|E(H)|=n-\binom{n}{k}$.  

To see that the difference between the burning number $b(H)$ and the lazy burning number $b_L(H)$ can be arbitrarily large (Corollary~\ref{large_differences}), we will consider {\em tight paths}.  Specifically, the $k$-uniform tight path of order $n$ is the hypergraph $H$ with $V(H)=\{v_1, \ldots, v_n\}$ and 
\[
E(H)=\big{\{} \{v_1,\dots,v_k\},\{v_2,\ldots,v_{k+1}\},\ldots,\{v_{n-k+1},\ldots,v_n\} \big{\}}.
\]

\begin{definition}
Given a $k$-uniform tight path $H$ on $n$ vertices, define a \emph{seed} as a set of $k-1$ vertices all belonging to a common edge. Define a \emph{burned seed} as a seed whose $k-1$ vertices were chosen as sources in the arsonist's burning sequence. 
\end{definition}

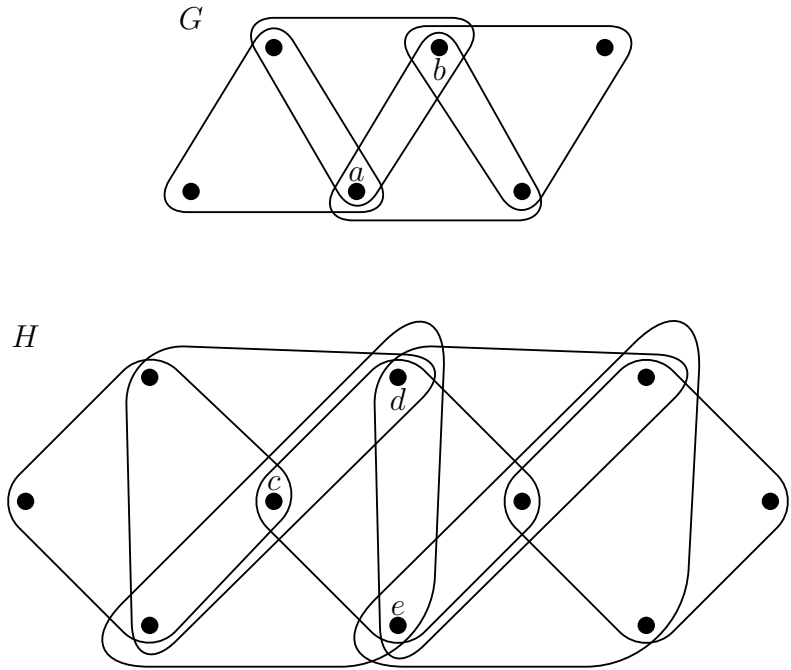
\begin{figure}[h]
\centering
\centering
\begin{tikzpicture}[scale=1.1]

\node () at (-2,2.1) {$G$};

\draw [line width=0.25mm] [black] [rounded corners=0.5cm] (-2.5,-0.25)--(-1,2.2)--(0.5,-0.25)--cycle;
\draw [line width=0.25mm] [black] [rounded corners=0.5cm] (-1.45,2.1)--(0,-0.4)--(1.6,2.1)--cycle;
\draw [line width=0.25mm] [black] [rounded corners=0.5cm] (-0.5,-0.35)--(1,2.15)--(2.4,-0.35)--cycle;
\draw [line width=0.25mm] [black] [rounded corners=0.5cm] (0.4,2)--(2,-0.45)--(3.5,2)--cycle;

\node (a) at (0,0) {};
\fill [fill=black] (a) circle (0.105) node [above] {$a$};
\node (a) at (-1,1.74) {};
\fill [fill=black] (a) circle (0.105) node [below] {};
\node (a) at (-2,0) {};
\fill [fill=black] (a) circle (0.105) node [above] {};
\node (a) at (1,1.74) {};
\fill [fill=black] (a) circle (0.105) node [below] {$b$};
\node (a) at (2,0) {};
\fill [fill=black] (a) circle (0.105) node [above] {};
\node (a) at (3,1.74) {};
\fill [fill=black] (a) circle (0.105) node [left] {};
\end{tikzpicture}

\begin{tikzpicture}[scale=1.1]

\node () at (-4.5,2) {$H$};
\node (a) at (0,3) {};
\fill [fill=white] (a) circle (0.105) node [above] {};

\draw [line width=0.25mm] [black] [rounded corners=0.5cm] (-4.9,0)--(-3,1.9)--(-1.1,0.1)--(-3,-1.9)--cycle;
\draw [line width=0.25mm] [black] [rounded corners=0.5cm] (-1.9,0)--(0,1.9)--(1.9,0)--(0,-1.9)--cycle;
\draw [line width=0.25mm] [black] [rounded corners=0.5cm] (1.1,0)--(3,1.9)--(4.9,0)--(3,-1.9)--cycle;

\draw [line width=0.25mm] [black] [rounded corners=0.8cm] (-3.2,-2.2)--(-3.3,1.9)--(0.8,1.75)--cycle;
\draw [line width=0.25mm] [black] [rounded corners=0.8cm] (-0.2,-2.25)--(-0.3,1.9)--(3.85,1.75)--cycle;

\draw [line width=0.25mm] [black] [rounded corners=1.2cm] (-4.1,-2)--(0.6,2.7)--(0.4,-2)--cycle;
\draw [line width=0.25mm] [black] [rounded corners=1.3cm] (-1.1,-2)--(3.7,2.75)--(3.45,-2)--cycle;

\node (a) at (-4.5,0) {};
\fill [fill=black] (a) circle (0.105) node [above] {};
\node (a) at (-1.5,0) {};
\fill [fill=black] (a) circle (0.105) node [above] {$c$};
\node (a) at (1.5,0) {};
\fill [fill=black] (a) circle (0.105) node [above] {};
\node (a) at (4.5,0) {};
\fill [fill=black] (a) circle (0.105) node [above] {};
\node (a) at (-3,1.5) {};
\fill [fill=black] (a) circle (0.105) node [above] {};
\node (a) at (0,1.5) {};
\fill [fill=black] (a) circle (0.105) node [below] {$d$};
\node (a) at (3,1.5) {};
\fill [fill=black] (a) circle (0.105) node [above] {};
\node (a) at (-3,-1.5) {};
\fill [fill=black] (a) circle (0.105) node [above] {};
\node (a) at (0,-1.5) {};
\fill [fill=black] (a) circle (0.105) node [above] {$e$};
\node (a) at (3,-1.5) {};
\fill [fill=black] (a) circle (0.105) node [above] {};
\end{tikzpicture}

\caption{A tight $3$-uniform path $G$ and a tight $4$-uniform path $H$.}
\label{T6}
\end{figure}

See Figure \ref{T6}. An example of a seed in $G$ is $\{a,b\}$, and an example of a seed in $H$ is  $\{c,d,e\}$.



Given a hypergraph $G$, define $S_G(r)$ as the maximum number of vertices that could possibly be on fire in $G$ at the end of round $r$ over all possible burning sequences. Note that the definition of $S_G(r)$ only applies to the round-based model of hypergraph burning.

\begin{lemma}
\label{obvious_lemma}
Let $G$ be a hypergraph and $r$ be the smallest natural number satisfying\linebreak 
$S_G(r)=|V(G)|$.
Then $r\leq b(G)$.
\end{lemma}

\begin{proof}
 Since $S_G(\ell)<|V(G)|$ for all $\ell<r$, $G$ cannot be fully burned in fewer than $r$ rounds.
 \end{proof}

\begin{theorem}
\label{T_n_lem}
For any tight $3$-uniform path $H$ on $n$ vertices, $b(H)=\left\lceil \sqrt{2n-1}\ \right\rceil$.
\end{theorem}


\begin{proof}
Write $V(H)=\{u_1,\ldots,u_n\}$. Clearly, no fire will propagate in $H$ without the existence of a burned seed. Also, once a burned seed $\{u_i, u_{i+1}\}$ is created in round $r$, it will cause fire to spread to $u_{i-1}$ and $u_{i+2}$ in round $r+1$ (if these vertices are not already on fire). In round $r+2$, it will cause fire to spread to $u_{i-2}$ and $u_{i+3}$, etc. The arsonist's best strategy is to create a burned seed every second round, spaced such that no vertex $u_j$ will catch fire due to propagation from two burned seeds until possibly the final round when $H$ is completely burned. By following this strategy, the maximum number of vertices that can catch fire through propagation  in round $r$ is twice the number of burned seeds that existed at the end of round $r-1$. Define $n(r)$ as the number of vertices that become burned in round $r$ when following this strategy (including sources), so $n(r)=S_H(r)-S_H(r-1)$. Clearly $n(1)=n(2)=1$ since no propagation can occur in the first two rounds. There is \emph{one} burned seed at the start of rounds three and four, so $n(3)=n(4)=2(1)+1=3$. Similarly, there are \emph{two} burned seeds at the start of rounds five and six, so $n(5)=n(6)=2(2)+1=5$. The sequence of $n(r)$-values is therefore $1,1,3,3,5,5,7,7,\ldots$ . Observe that $S_H(r)=\sum_{i=1}^r n(i)$, so the sequence of $S_H(r)$-values is $1,2,5,8,13,18,25,32,41,50,\ldots$ . The formula for this sequence is known to be $S_H(r)=\left\lfloor \frac{r^2+1}{2} \right\rfloor$; see sequence A000982 in \cite{oeis}. 

Let $x_n$ be the smallest natural number satisfying $S_H(x_n) \geq |V(H)|=n$, so $x_n\leq b(H)$ by Lemma \ref{obvious_lemma}. In the case of our hypergraph $H$, $b(H)$ is exactly equal to $x_n$, since the arsonist can successfully burn $H$ in $x_n$ rounds using the strategy outlined above. We will now show that $x_n=\left\lceil \sqrt{2n-1}\ \right\rceil$. 

Suppose $x_n<\sqrt{2n-1}$. Then we have $$S_H(x_n)=\left\lfloor \frac{x_n^2+1}{2} \right\rfloor \leq \frac{x_n^2+1}{2} < \frac{(\sqrt{2n-1})^2+1}{2}=n,$$ which is a contradiction. We therefore have the restriction $x_n\geq \sqrt{2n-1}$, and the least integer that fulfills this restriction is $\left\lceil \sqrt{2n-1}\ \right\rceil$. Indeed, since $\left\lceil \sqrt{2n-1}\ \right\rceil\geq \sqrt{2n-1}$ and $S_H(r)=\left\lfloor \frac{r^2+1}{2} \right\rfloor$ is a monotonic increasing function, we have $$S_H(\left\lceil \sqrt{2n-1}\ \right\rceil)=\left\lfloor \frac{\left\lceil \sqrt{2n-1}\ \right\rceil^2+1}{2} \right\rfloor \geq \left\lfloor \frac{(\sqrt{2n-1})^2+1}{2} \right\rfloor =\lfloor n\rfloor=n.$$ Therefore $x_n=\left\lceil \sqrt{2n-1}\ \right\rceil$, and hence $b(H)=\left\lceil \sqrt{2n-1}\ \right\rceil$.
\end{proof}

\begin{corollary}
\label{large_differences}
Given any $k\in\mathbb{N}$, there exist hypergraphs $G$ and $H$ such that $b(G)-b_L(G)>k$ and $\frac{b(H)}{b_L(H)}>k$.
\end{corollary}

\begin{proof}
Denote by $H_n$ the $3$-regular tight path on $n$ vertices. For all $n\geq 3$, $b_L(H_n)=2$. But $b(H_n)=\left\lceil \sqrt{2n-1}\ \right\rceil$. Thus, $b(H_n)-b_L(H_n)=\left\lceil \sqrt{2n-1}\ \right\rceil-2$ and $\frac{b(H_n)}{b_L(H_n)}=\frac{\left\lceil \sqrt{2n-1}\ \right\rceil}{2}$, which are both unbounded.
\end{proof}

We will describe the optimal burning sequence for a tight 3-uniform path $H_n$ on $n$ vertices. First, we recall how to find the optimal burning sequence on a path $P_n$ from \cite{first_paper}, as this is the inspiration for how we burn $H_n$. In particular, $b(P_n)=\left\lceil\sqrt{n}\ \right\rceil$ \cite{first_paper}, and if $n$ is a square number, $P_n$ can be burned ``perfectly'' (i.e.~with no overlap of propagating fires) in $\sqrt{n}$ rounds as follows. Write $V(P_n)=\{u_1,\ldots,u_n\}$ and $E(P_n)=\big{\{} \{u_1,u_2\},\ldots,\{u_{n-1},u_n\}\big{\}}$. Set $k=\sqrt{n}$, and denote the burning sequence by $(x_1,\ldots,x_k)$. Finally, for $i=0,1,\ldots,k-1$, set $x_{k-i}=u_{n-i^2-i}$. If $n$ is not a square number, let $N$ be the smallest square number exceeding $n$, so $b(P_n)=b(P_N)$. Consider an optimal burning sequence on $P_N$ as described above. Then, delete the $N-n$ lowest-indexed vertices $u_i$ to create $P_n$ from $P_N$. If the first source $x_1$ was deleted, the new first source becomes $u_{1+N-n}$ (i.e.~the new lowest-indexed vertex in the resulting path).

Our strategy for burning $H_n$ is similar. If $n$ is twice a square number \emph{or} $\sqrt{2n-1}\in\mathbb{N}$, then $H_n$ can be burned ``perfectly'' in $\left\lceil\sqrt{2n-1}\ \right\rceil$ rounds as follows. Write $V(H_n)=\{u_1,\ldots,u_n\}$ and $E(H_n)=\big{\{} \{u_1,u_2,u_3\},\ldots,\{u_{n-2},u_{n-1},u_n\}\big{\}}$. Set $k=\left\lceil\sqrt{2n-1}\ \right\rceil$ and denote the burning sequence by $(x_1,\ldots,x_k)$. Note that if $n$ is twice a square number then $k$ will be even, and if $\sqrt{2n-1}\in\mathbb{N}$ then $k$ will be odd. Now, for $s=1,2,\ldots,\left\lceil\frac{k}{2}\right\rceil$, the $s^{th}$ odd-indexed source is $x_{2s-1}=u_a$, where $$a=k-2s+2+2\sum_{i=1}^{s-1}(k-2i+1)\ .$$
Similarly, for $s=1,2,\ldots,\left\lfloor\frac{k}{2}\right\rfloor$, the $s^{th}$ even-indexed source is $x_{2s}=u_b$, where $$b=k-2s+1+2\sum_{i=1}^{s-1}(k-2i+1)\ .$$ Of course, $x_{2s-1}$ and $x_{2s}$ make up the $s^{th}$ seed. These formulas were found using a technique similar to the one used in \cite{first_paper}. In particular, the $s^{th}$ seed causes all vertices within distance $k-2s$ to catch fire by the end of round $k$. Spacing the seeds ``perfectly'' so that there is no overlap in propagating fires gives the above formulas, as well as the fact that $n$ is either twice a square number \emph{or} $\sqrt{2n-1}\in\mathbb{N}$. Finally, if $n$ has neither of these properties, then let $N$ be the smallest number exceeding $n$ that is twice a square number \emph{or} satisfies $\sqrt{2N-1}\in\mathbb{N}$. Then $b(H_n)=b(H_N)$. Consider an optimal burning sequence on $H_N$ as described above, then delete the $N-n$ lowest-indexed vertices $u_i$ to create $H_n$ from $H_N$. If $x_1$ and/or $x_2$ were deleted, then the new first seed becomes $u_{1+N-n}$ and $u_{2+N-n}$.

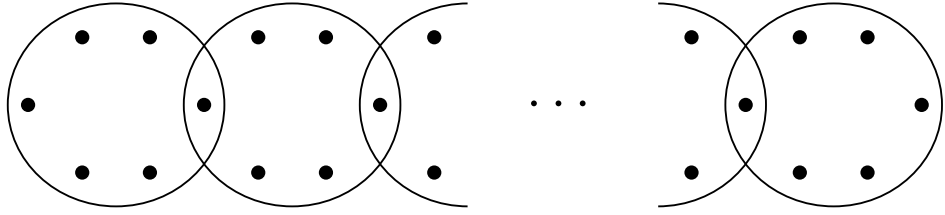
\begin{figure}[h]
\centering
\begin{tikzpicture}[scale=0.9]

\node (d) at (0,0) {};
\fill [fill=black] (d) circle (0.105) node [above] {};
\node (b) at (-0.8,1) {};
\fill [fill=black] (b) circle (0.105) node [below left] {};   
\node (e) at (-1.8,1) {};
\fill [fill=black] (e) circle (0.105) node [left] {};    
\node (a) at (-2.6,0) {};
\fill [fill=black] (a) circle (0.105) node [below right] {};
\node (c) at (-1.8,-1) {};
\fill [fill=black] (c) circle (0.105) node [right] {};
\node (f) at (-0.8,-1) {};
\fill [fill=black] (f) circle (0.105) node [right] {};   

\node (b) at (0.8,1) {};
\fill [fill=black] (b) circle (0.105) node [below left] {};   
\node (b) at (0.8,-1) {};
\fill [fill=black] (b) circle (0.105) node [below left] {};   

\node (b) at (-3.4,1) {};
\fill [fill=black] (b) circle (0.105) node [below left] {};   
\node (e) at (-4.4,1) {};
\fill [fill=black] (e) circle (0.105) node [left] {};    
\node (a) at (-5.2,0) {};
\fill [fill=black] (a) circle (0.105) node [below right] {};
\node (c) at (-4.4,-1) {};
\fill [fill=black] (c) circle (0.105) node [right] {};
\node (f) at (-3.4,-1) {};
\fill [fill=black] (f) circle (0.105) node [right] {};   

\draw [line width=0.25mm, black ] (-1.3,0) ellipse (1.6cm and 1.5cm);
\draw [line width=0.25mm, black ] (-3.9,0) ellipse (1.6cm and 1.5cm);
\draw [line width=0.25mm, black ] (1.3,0) ellipse (1.6cm and 1.5cm);

\node (d) at (8,0) {};
\fill [fill=black] (d) circle (0.105) node [above] {};
\node (b) at (7.2,1) {};
\fill [fill=black] (b) circle (0.105) node [below left] {};   
\node (e) at (6.2,1) {};
\fill [fill=black] (e) circle (0.105) node [left] {};    
\node (a) at (5.4,0) {};
\fill [fill=black] (a) circle (0.105) node [below right] {};
\node (c) at (6.2,-1) {};
\fill [fill=black] (c) circle (0.105) node [right] {};
\node (f) at (7.2,-1) {};
\fill [fill=black] (f) circle (0.105) node [right] {};   

\node (a) at (4.6,1) {};
\fill [fill=black] (a) circle (0.105) node [below right] {};
\node (a) at (4.6,-1) {};
\fill [fill=black] (a) circle (0.105) node [below right] {};

\draw [line width=0.25mm, black ] (6.7,0) ellipse (1.6cm and 1.5cm);
\draw [line width=0.25mm, black ] (4.1,0) ellipse (1.6cm and 1.5cm);

\draw [fill=white,white] (1.3,1.8) rectangle (4.1,-1.8); 

\node (b) at (2.7,0) {\LARGE{\textbf{$\cdots$}}};
\end{tikzpicture}

\caption{A family of uniform, linear hypergraphs $H$ for which no upper bounds exist on $b(H)$ or $b_L(H)$ that are sublinear in terms of $|V(H)|$.}
\label{ex_path}
\end{figure}

Given that hypergraph burning is an extension of graph burning, it is natural to ask if there is an analogous conjecture for hypergraphs to the Burning Number Conjecture for graphs. That is, for arbitrary hypergraphs $H$, we wish to find a sublinear upper bound on $b(H)$ (and $b_L(H)$) in terms of $|V(H)|$. It turns out that no such bound exists even when we insist $H$ is $k$-uniform and linear. Consider the family of $k$-uniform loose paths, which are indeed linear; see Figure \ref{ex_path}. The best possible upper bounds for this family of hypergraphs are $b_L(H),b(H) \in \mathcal{O}(|V(H)|)$. To see this, recall that $|V(H)|-|E(H)|\leq b_L(H)<b(H)\leq|V(H)|$, and observe that we can choose the size of the edges large enough so that $|V(H)|-|E(H)|\in\Theta(|V(H)|)$. 



Furthermore, for any \emph{fixed} $k$, this family of hypergraphs has no sublinear upper bound on its (lazy) burning number in terms of its order. Indeed, any $k$-uniform loose path $H$ has $|E(H)|=\frac{|V(H)|-1}{k-1}$. Thus, a \emph{lower} bound on $b_L(H)$ and $b(H)$ is $$|V(H)|-|E(H)|=|V(H)|-\left(\frac{|V(H)|-1}{k-1}\right)=\left(\frac{k-2}{k-1}\right)|V(H)|+\frac{1}{k-1}\ ,$$ which is linear in $|V(H)|$.


We close this section with a brief discussion of complexity. It is easy to come up with a polynomial-time algorithm which takes a sequence of vertices in a hypergraph $H$ as input, and determines if it is a valid burning sequence for $H$. There is a similar algorithm for the lazy burning game which determines if a given subset of $V(H)$ is a lazy burning set in polynomial time. Hence, both games are in \textbf{NP}. Since graph burning is \textbf{NP}-complete \cite{burning_is_hard}, and it is a special case of hypergraph burning, hypergraph burning is also \textbf{NP}-complete. We leave the question of whether lazy hypergraph burning is \textbf{NP}-complete as an open problem. 

\section{Disconnected Hypergraphs}
\label{disconnected_section}

In this section we consider how to write the 
burning number of a disconnected hypergraph in terms of the 
burning numbers of its connected components. The solution is trivial in the lazy case (see Lemma \ref{H_components_lazy}), but the round-based case is more complicated. 

\begin{lemma}
\label{H_components_lazy}
If $H$ is disconnected with connected components $G_1,G_2,\ldots,G_k$, then $b_L(H)=b_L(G_1)+b_L(G_2)+ \cdot\cdot\cdot +b_L(G_k)$.
\end{lemma}

\begin{proof}
Clearly if $S_i$ is a minimum lazy burning set for $G_i$ for each $i\in\{1,2,\ldots,k\}$, then a minimum lazy burning set for $H$ is $S_1\cup S_2\cup \cdot\cdot\cdot \cup S_k$.
\end{proof}

\begin{figure}[h]
\centering
\hfill
\begin{subfigure}{0.45\textwidth}
\centering
\begin{tikzpicture}[scale=1.1]

\draw [line width=0.25mm] [black] [rounded corners=0.5cm] (-2.5,-0.25)--(-1,2.2)--(0.5,-0.25)--cycle;
\draw [line width=0.25mm] [black] [rounded corners=0.5cm] (-1.45,2.1)--(0,-0.4)--(1.6,2.1)--cycle;
\draw [line width=0.25mm] [black] [rounded corners=0.5cm] (-0.6,-0.35)--(1,2.15)--(2.4,-0.35)--cycle;
\draw [line width=0.25mm] [black] [rounded corners=0.5cm] (0.4,2)--(2,-0.45)--(3.5,2)--cycle;

\node (a) at (0,0) {};
\fill [fill=black] (a) circle (0.105) node [above] {3};
\node (a) at (-1,1.74) {};
\fill [fill=black] (a) circle (0.105) node [below right] {2};
\node (a) at (-2,0) {};
\fill [fill=black] (a) circle (0.105) node [right] {1};
\node (a) at (1,1.74) {};
\fill [fill=black] (a) circle (0.105) node [below] {4};
\node (a) at (2,0) {};
\fill [fill=black] (a) circle (0.105) node [above left] {5};
\node (a) at (3,1.74) {};
\fill [fill=black] (a) circle (0.105) node [left] {6};
\end{tikzpicture}
\end{subfigure}
\begin{subfigure}{0.45\textwidth}
\begin{tikzpicture}[scale=1.1]

\draw [line width=0.25mm] [black] [rounded corners=0.5cm] (-2.5,-0.25)--(-1,2.2)--(0.5,-0.25)--cycle;
\draw [line width=0.25mm] [black] [rounded corners=0.5cm] (-1.45,2.1)--(0,-0.4)--(1.65,2)--cycle;
\draw [line width=0.25mm] [black] [rounded corners=0.5cm] (-0.6,-0.35)--(1,2.15)--(2.4,-0.3)--cycle;

\node (a) at (0,0) {};
\fill [fill=black] (a) circle (0.105) node [above] {9};
\node (a) at (-1,1.74) {};
\fill [fill=black] (a) circle (0.105) node [below right] {8};
\node (a) at (-2,0) {};
\fill [fill=black] (a) circle (0.105) node [right] {7};
\node (a) at (1,1.74) {};
\fill [fill=black] (a) circle (0.105) node [below] {10};
\node (a) at (2,0) {};
\fill [fill=black] (a) circle (0.105) node [left] {11};
\end{tikzpicture}
\end{subfigure}

\caption{A hypergraph whose burning number is strictly less than the sum of the burning numbers of its connected components.}
\label{strict_ineq_ex}
\end{figure}
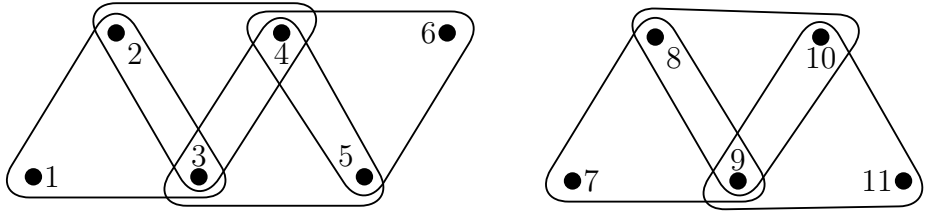

We cannot say in general that a similar equality to the one in Lemma \ref{H_components_lazy} holds for $b(H)$. It is indeed possible that $b(H)<b(G_1)+b(G_2)+ \cdot\cdot\cdot +b(G_k)$. The arsonist may ``save time'' by initially burning enough vertices to start propagation in one component, then moving on to subsequent components. For example, consider the disconnected hypergraph $H$ in Figure~\ref{strict_ineq_ex}. The ``left'' component of $H$ has burning number $4$, and the ``right'' component of $H$ has burning number $3$. However, $(3,4,9,10,7)$ is an optimal burning sequence, so $b(H)=5$, which is less than the sum of the burning numbers of the two connected components.


\begin{lemma}
\label{easy_discon_lem_for_b}
If $H$ is disconnected with connected components $G_1,G_2,\ldots,G_k$, then\\ $\max\{b(G_1), b(G_2),\ldots,b(G_k)\}\leq b(H)\leq b(G_1)+b(G_2)+ \cdot\cdot\cdot +b(G_k)$.
\end{lemma}

The proof of Lemma \ref{easy_discon_lem_for_b} is omitted, as both inequalities are intuitive (see \cite{mythesis} for a proof). We can improve the upper bound from Lemma \ref{easy_discon_lem_for_b} if we assume none of the connected components of $H$ are isolated vertices.

\begin{lemma}
\label{sum_with_minus_k}
If $H$ is disconnected with connected components $G_1,G_2,\ldots,G_k$, none of which are isolated vertices, then $b(H)\leq b(G_1)+b(G_2)+ \cdot\cdot\cdot +b(G_k)-k+1$.
\end{lemma}

\begin{proof}
For each $G_i$ let $S_i=(u^{(i)}_1,\ldots,u^{(i)}_{b(G_i)})$ be an optimal burning sequence. We claim $$S=(u^{(1)}_1,\ldots,u^{(1)}_{b(G_1)-1},u^{(2)}_1,\ldots,u^{(2)}_{b(G_2)-1},\ldots,u^{(k)}_1,\ldots,u^{(k)}_{b(G_k)-1},u^{(k)}_{b(G_k)})$$ is a burning sequence for $H$. For each $i \in \{1,2,\ldots,k-1\}$ burning $S_i \setminus \{u^{(i)}_{b(G_i)}\}$ one-by-one and in order will leave the entirety of $G_i$ burned in at most two more rounds via propagation (to see this, consider the two possible cases: $u^{(i)}_{b(G_i)}$ is a redundant source or a non-redundant source). Since $G_k$ is not an isolated vertex, it has burning number at least two, so there is ``enough time'' for $G_{k-1}$ to fully catch fire while the arsonist burns $G_k$. Since $S$ is a valid burning sequence for $H$ we have $b(H)\leq |S|=(b(G_1)-1)+\cdot\cdot\cdot + (b(G_{k-1})-1)+b(G_k)=b(G_1)+\cdot\cdot\cdot +b(G_k)-k+1$.
\end{proof}

The upper bound from Lemma \ref{sum_with_minus_k} is tight. Consider the hypergraph $H$ with $V(H)=\{1,2,.\ldots,21\}$ and $E(H)=\big{\{}\{1,2,3\},\{4,5,6\},\ldots,\{19,20,21\}\big{\}}$. The subhypergraph induced by \emph{one} of the seven edges has burning number $3$, and since $(1,2,4,5,7,8,10,11,13,14,$ $16,17,19,20,21)$ is an optimal burning sequence for $H$, it has burning number $15=(3)(7)-7+1$. Indeed, this example can be expanded to an infinite family of hypergraphs which exhibit the tightness of the bound in Lemma \ref{sum_with_minus_k}. Such a hypergraph $G$ consists of multiple disjoint non-singleton edges $e_1,\ldots,e_k$, with no other edges and no isolated vertices. For each $i\in\{1,\ldots,k\}$, let $G_i$ be the subhypergraph of $G$ induced by $e_i$. Then, one may construct an optimal burning sequence for $G$ by burning $|e_i|-1=b(G_i)-1$ sources in $e_i$ for each $i\in\{1,\ldots,k-1\}$, followed by $|e_k|=b(G_k)$ sources in $e_k$ (where the last source is redundant). This shows the tightness of the bound in Lemma  \ref{sum_with_minus_k}, since the total number of sources is $$\big{(}b(G_1)-1\big{)}+\cdots+\big{(}b(G_{k-1})-1\big{)}+b(G_k)= b(G_1)+ \cdots +b(G_k)-k+1 .$$

\begin{figure}[h]
\centering

\begin{tikzpicture}[scale=1.1]

\draw [line width=0.25mm] [black] [rounded corners=0.5cm] (-2.5,-0.25)--(-1,2.2)--(0.5,-0.25)--cycle;
\draw [line width=0.25mm] [black] [rounded corners=0.5cm] (-1.45,2.1)--(0,-0.4)--(1.6,2.1)--cycle;
\draw [line width=0.25mm] [black] [rounded corners=0.5cm] (-0.6,-0.35)--(1,2.15)--(2.6,-0.4)--cycle;
\draw [line width=0.25mm] [black] [rounded corners=0.5cm] (0.4,2)--(2,-0.45)--(3.65,2.05)--cycle;

\node (a) at (0,0) {};
\fill [fill=black] (a) circle (0.105) node [above] {};
\node (a) at (-1,1.74) {};
\fill [fill=black] (a) circle (0.105) node [below right] {};
\node (a) at (-2,0) {};
\fill [fill=black] (a) circle (0.105) node [right] {};
\node (a) at (1,1.74) {};
\fill [fill=black] (a) circle (0.105) node [below] {};
\node (a) at (2,0) {};
\fill [fill=black] (a) circle (0.105) node [above] {$u_1$};
\node (a) at (3,1.74) {};
\fill [fill=black] (a) circle (0.105) node [below] {$u_2$};
\node (a) at (4,0) {};
\fill [fill=black] (a) circle (0.105) node [above] {};
\node (a) at (5,1.74) {};
\fill [fill=black] (a) circle (0.105) node [below] {};
\node (a) at (6,0) {};
\fill [fill=black] (a) circle (0.105) node [above] {};
\node (a) at (7,1.74) {};
\fill [fill=black] (a) circle (0.105) node [below] {};

\draw [line width=0.25mm] [black] [rounded corners=0.5cm] (1.35,-0.3)--(3,2.2)--(4.5,-0.25)--cycle;
\draw [line width=0.25mm] [black] [rounded corners=0.5cm] (2.35,2.15)--(4,-0.4)--(5.6,2.1)--cycle;
\draw [line width=0.25mm] [black] [rounded corners=0.5cm] (3.4,-0.35)--(5,2.15)--(6.4,-0.35)--cycle;
\draw [line width=0.25mm] [black] [rounded corners=0.5cm] (4.4,2)--(6,-0.45)--(7.5,2)--cycle;

\draw [line width=0.25mm] [black] [rounded corners=0.5cm] (7.5,-0.25)--(9,2.2)--(10.5,-0.25)--cycle;

\node (a) at (8,0) {};
\fill [fill=black] (a) circle (0.105) node [above] {};
\node (a) at (9,1.74) {};
\fill [fill=black] (a) circle (0.105) node [below] {$u_3$};
\node (a) at (10,0) {};
\fill [fill=black] (a) circle (0.105) node [left] {$u_4$};

\node (a) at (11.5,1.74) {};
\fill [fill=black] (a) circle (0.105) node [below] {$u_5$};
\node (a) at (11.5,0) {};
\fill [fill=black] (a) circle (0.105) node [above] {$u_6$};

\draw [line width=0.25mm, black ] (11.5,0.87) ellipse (0.4cm and 1.15cm);

\end{tikzpicture}

\caption{An example where the inequality in Lemma \ref{sum_with_minus_k} is strict.}
\label{disconn3}
\end{figure}
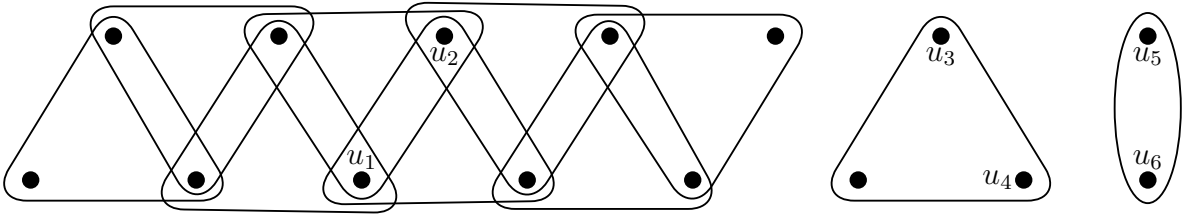

There are also examples where the upper bound from Lemma \ref{sum_with_minus_k} is a strict inequality. Consider the hypergraph in Figure \ref{disconn3}. The burning numbers of the left, middle, and right connected components are $5$, $3$, and $2$ respectively. But $(u_1,u_2,u_3,u_4,u_5,u_6)$ is an optimal burning sequence for the entire hypergraph, so it has burning number $6<(5+3+2)-3+1$. It is an open question as to whether or not there is an infinite family of hypegraphs for which the bound in Lemma \ref{sum_with_minus_k} is strict.



\section{Subhypergraphs}
\label{subhypergraph_section}

In this section we compare the (lazy) burning number of a hypergraph $H$ to the (lazy) burning numbers of various types of subhypergraphs of $H$. Before describing the subhypergraphs of interest, we first define the \emph{incidence matrix} of a hypergraph $H=(V,E)$, which is the $|V|\times |E|$ matrix whose $(i,j)$-entry is a 1 if the $i^{th}$ vertex belongs to the $j^{th}$ edge, and a 0 otherwise.


A hypergraph $H^\prime=(V^\prime,E^\prime)$ is called a \emph{weak subhypergraph} (or just \emph{subhypergraph}) of \linebreak$H=(V,E)$ if $V^\prime \subseteq V$ and either $E^\prime=\emptyset$, or, after a suitable permutation of its rows and columns, the incidence matrix for $H^\prime$ is a submatrix of the incidence matrix of $H$. 
Thus, each edge $e^\prime \in E^\prime$ has $e^\prime=e \cap V^\prime$ for some $e\in E$ (so $H^\prime$ may not be simple, even if $H$ is). If, in addition, $E^\prime = \{e \cap V^\prime\mid e\in E, e \cap V^\prime\neq \emptyset \}$, then $H^\prime$ is said to be \emph{induced} by $V^\prime$ (that is, all nonempty edges $e \cap V^\prime$, singletons included, are present in $E^\prime$). 

A hypergraph $H^{\prime\prime}=(V^{\prime\prime},E^{\prime\prime})$ is called a \emph{strong subhypergraph} (or \emph{hypersubgraph}) of \linebreak$H=(V,E)$ if $V^{\prime\prime} \subseteq V$ and $E^{\prime\prime} \subseteq E$. If all edges $e\in E$ such that $e \subseteq V^{\prime\prime}$ 
are present in $E^{\prime\prime}$, then $H^{\prime\prime}$ is said to be \emph{induced} by $V^{\prime\prime}$, and we write $H^{\prime\prime}=H[V^{\prime\prime}]$. If $V^{\prime\prime} = \cup_{e \in E^{\prime\prime}} e$ then $H^{\prime\prime}$ is said to be \emph{induced} by $E^{\prime\prime}$, and we write $H^{\prime\prime}=H[E^{\prime\prime}]$. Note that every strong subhypergraph of $H$ is also a weak subhypergraph of $H$, but the converse is not necessarily true.

Note that we allow the existence  of singleton edges, as these often arise when taking a weak induced subhypergraph. Several results in this section use a ``shadow strategy'' proof technique which involves comparing two different instances of hypergraph burning. For the sake of avoiding confusion during these proofs, we will denote the burning game that takes place on a hypergraph $H$ when burning according to the sequence $S$ by $BG(H,S)$. Similarly,  we will denote the lazy burning game that takes place on a hypergraph $H$ with lazy burning set $L$ by $LBG(H,L)$.

\begin{lemma}
\label{same_vertex_set}
Let $G_1$ and $G_2$ be two hypergraphs with the same vertex set such that $E(G_1)\subseteq E(G_2)$. Then $b_L(G_2)\leq b_L(G_1)$ and $b(G_2)\leq b(G_1)$.
\end{lemma}



The proof of Lemma \ref{same_vertex_set} is omitted since it is quite intuitive; $G_2$ can be created from $G_1$ by adding edges, so $G_2$ is ``more flammable'' than $G_1$ (see \cite{mythesis} for a full proof). Our next result compares the (lazy) burning numbers of a weak and a strong subhypergraph induced by the same set of vertices. 

\begin{lemma}
\label{strvswk}
Let $H$ be a hypergraph with $V(H)=\{v_1,\ldots,v_n\}$ and let $I\subseteq \{1,\ldots,n\}$. If $G_1$ is a strong subhypergraph induced by $\{v_i \mid i\in I\}$ and $G_2$ is a weak subhypergraph induced by $\{v_i \mid i\in I\}$ then $b(G_2)\leq b(G_1)$ and $b_L(G_2)\leq b_L(G_1)$.
\end{lemma}

\begin{proof}
Both inequalities follow from Lemma \ref{same_vertex_set} and the fact that $V(G_1)=V(G_2)$ and $E(G_1)\subseteq E(G_2)$.
\end{proof}

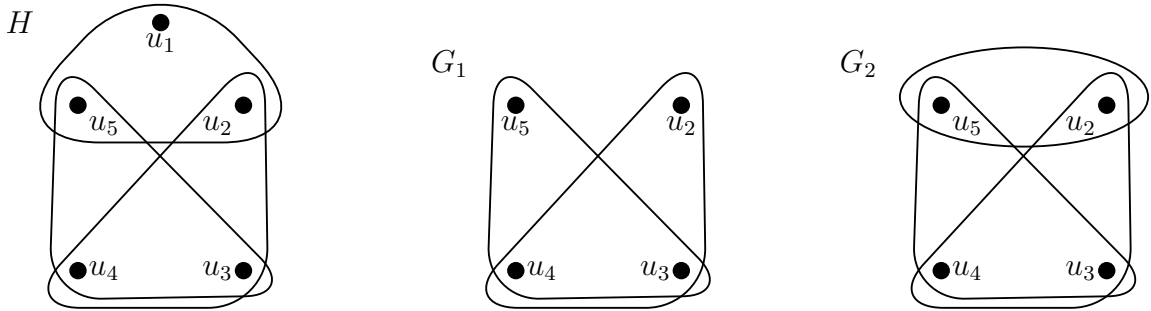
\begin{figure}[h]
\centering
\hfill
\begin{subfigure}{0.3\textwidth}
\centering
\begin{tikzpicture}[scale=1.1]

\node () at (-1.7,2) {$H$};

\draw [line width=0.25mm] [black] [rounded corners=0.7cm] (-1.25,1.65)--(1.65,-1.3)--(-1.35,-1.35)--cycle;
\draw [line width=0.25mm] [black] [rounded corners=0.8cm] (1.25,1.75)--(-1.7,-1.45)--(1.3,-1.45)--cycle;
\draw [line width=0.25mm] [black] [rounded corners=1.5cm] (-2.1,0.55)--(2.1,0.55)--(0,2.8)--cycle;

\node (a) at (0,2) {};
\fill [fill=black] (a) circle (0.105) node [below] {$u_1$};
\node (a) at (1,1) {};
\fill [fill=black] (a) circle (0.105) node [below left] {$u_2$};
\node (a) at (1,-1) {};
\fill [fill=black] (a) circle (0.105) node [left] {$u_3$};
\node (a) at (-1,-1) {};
\fill [fill=black] (a) circle (0.105) node [right] {$u_4$};
\node (a) at (-1,1) {};
\fill [fill=black] (a) circle (0.105) node [below right] {$u_5$};
\end{tikzpicture}

\end{subfigure}
\hfill
\begin{subfigure}{0.3\textwidth}
\centering
\begin{tikzpicture}[scale=1.1]

\node () at (-1.8,1.5) {$G_1$};
\node (a) at (0,2.5) {};
\fill [fill=white] (a) circle (0.105) node [below] {};

\draw [line width=0.25mm] [black] [rounded corners=0.7cm] (-1.25,1.65)--(1.65,-1.3)--(-1.35,-1.35)--cycle;
\draw [line width=0.25mm] [black] [rounded corners=0.8cm] (1.25,1.75)--(-1.7,-1.45)--(1.3,-1.45)--cycle;

\node (a) at (1,1) {};
\fill [fill=black] (a) circle (0.105) node [below] {$u_2$};
\node (a) at (1,-1) {};
\fill [fill=black] (a) circle (0.105) node [left] {$u_3$};
\node (a) at (-1,-1) {};
\fill [fill=black] (a) circle (0.105) node [right] {$u_4$};
\node (a) at (-1,1) {};
\fill [fill=black] (a) circle (0.105) node [below] {$u_5$};
\end{tikzpicture}

\end{subfigure}
\hfill
\begin{subfigure}{0.3\textwidth}
\centering
\begin{tikzpicture}[scale=1.1]

\node () at (-2,1.5) {$G_2$};
\node (a) at (0,2.5) {};
\fill [fill=white] (a) circle (0.105) node [below] {};

\draw [line width=0.25mm] [black] [rounded corners=0.7cm] (-1.25,1.65)--(1.65,-1.3)--(-1.35,-1.35)--cycle;
\draw [line width=0.25mm] [black] [rounded corners=0.8cm] (1.25,1.75)--(-1.7,-1.45)--(1.3,-1.45)--cycle;

\node (a) at (1,1) {};
\fill [fill=black] (a) circle (0.105) node [below left] {$u_2$};
\node (a) at (1,-1) {};
\fill [fill=black] (a) circle (0.105) node [left] {$u_3$};
\node (a) at (-1,-1) {};
\fill [fill=black] (a) circle (0.105) node [right] {$u_4$};
\node (a) at (-1,1) {};
\fill [fill=black] (a) circle (0.105) node [below right] {$u_5$};

\draw [line width=0.25mm, black ] (0,1.1) ellipse (1.5cm and 0.6cm);
\end{tikzpicture}

\end{subfigure}

\caption{An example showing that equality can hold in Lemma \ref{strvswk}.}
\label{strvswk2}
\end{figure}

The bounds in Lemma \ref{strvswk} are both tight. Consider the hypergraph $H$ in Figure \ref{strvswk2}. The subset of vertices $\{u_2,u_3,u_4,u_5\}$ induces the strong subhypergraph $G_1$ and the weak subhypergraph $G_2$. Indeed, $\{u_2,u_3\}$ is a minimum lazy burning set in both $G_1$ and $G_2$. The subhypergraphs also have the same burning number, as $(u_2,u_3,u_5)$ is an optimal burning sequence in $G_1$, and $(u_2,u_3,u_4)$ is an optimal burning sequence $G_2$. Indeed, we may expand this to an infinite family of examples by adding any number of degree-one vertices to the edge containing $u_1$ in $H$ (\emph{without} renaming $u_1$ through $u_5$). Then, the strong and weak subhypergraphs induced by $\{u_2,u_3,u_4,u_5\}$ will be $G_1$ and $G_2$ (from Figure \ref{strvswk2}) respectively.


\begin{figure}[h]
\centering
\hfill
\begin{subfigure}{0.3\textwidth}
\centering
\begin{tikzpicture}[scale=1.1]

\node () at (-1.5,1.5) {$H$};

\node (x) at (0,1) {};
\fill [fill=black] (x) circle (0.105) node [below] {$u_2$};
\node (y) at (1,0) {};
\fill [fill=black] (y) circle (0.105) node [left] {$u_4$};
\node (z) at (0,-1) {};
\fill [fill=black] (z) circle (0.105) node [above] {$u_3$};
\node (w) at (-1,0) {};
\fill [fill=black] (w) circle (0.105) node [right] {$u_1$};

\draw [line width=0.25mm] [black] [rounded corners=0.7cm] (0.25,1.5)--(0.25,-1.5)--(-1.5,0)--cycle;
\draw [line width=0.25mm] [black] [rounded corners=1cm] (-0.35,1.8)--(-0.35,-1.8)--(1.7,0)--cycle;
\end{tikzpicture}

\end{subfigure}
\hfill
\begin{subfigure}{0.3\textwidth}
\centering
\begin{tikzpicture}[scale=1.1]

\node () at (-1.5,1.5) {$G_1$};

\node (x) at (0,1) {};
\fill [fill=black] (x) circle (0.105) node [below] {$u_2$};
\node (z) at (0,-1) {};
\fill [fill=black] (z) circle (0.105) node [above] {$u_3$};
\node (w) at (-1,0) {};
\fill [fill=black] (w) circle (0.105) node [right] {$u_1$};

\draw [line width=0.25mm] [black] [rounded corners=0.7cm] (0.25,1.5)--(0.25,-1.5)--(-1.5,0)--cycle;
\end{tikzpicture}

\end{subfigure}
\hfill
\begin{subfigure}{0.3\textwidth}
\centering
\begin{tikzpicture}[scale=1.1]

\node () at (-1.5,1.5) {$G_2$};

\node (x) at (0,1) {};
\fill [fill=black] (x) circle (0.105) node [below] {$u_2$};
\node (z) at (0,-1) {};
\fill [fill=black] (z) circle (0.105) node [above] {$u_3$};
\node (w) at (-1,0) {};
\fill [fill=black] (w) circle (0.105) node [right] {$u_1$};

\draw [line width=0.25mm] [black] [rounded corners=0.7cm] (0.25,1.5)--(0.25,-1.5)--(-1.5,0)--cycle;
\draw [line width=0.25mm, black ] (0.05,0) ellipse (0.4cm and 1.4cm);
\end{tikzpicture}

\end{subfigure}

\caption{An example showing that the inequality in Lemma \ref{strvswk} can be strict.}
\label{strvswk1}
\end{figure}
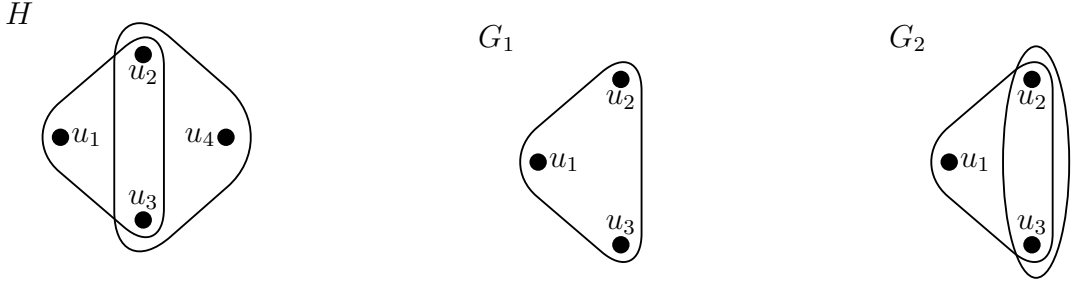

The inequalities in Lemma \ref{strvswk} can also be strict. Consider the hypergraph $H$ in Figure~\ref{strvswk1}. The subset of vertices $\{u_1,u_2,u_3\}$ induces the strong subhypergraph $G_1$ and the weak subhypergraph $G_2$. Indeed, $b_L(G_2)=1<2=b_L(G_1)$, and $b(G_2)=2<3=b(G_1)$. Again, we can expand this example to an infinite family of hypergraphs $H$ in which the inequality in Lemma \ref{strvswk} is strict. Construct such a hypergraph $H$ by letting $V(H)=\{u_1,\ldots,u_{k+1}\}$ and $E(H)=\big{\{}\{u_1,\ldots,u_k\},\{u_2,\ldots,u_{k+1}\}\big{\}}$ for some $k\geq3$. Let $G_1$ and $G_2$ be the strong and weak subhypergraphs respectively that are induced by $\{u_1,\ldots,u_k\}$. Then, $b(G_1)=k$ and $b_L(G_1)=k-1$ since $G_1$ consists only of $k$ vertices in an edge. However, $b(G_2)=k-1$ and $b_L(G_2)=k-2$. To see this, observe that $E(G_2)=\big{\{}\{u_1,\ldots,u_k\},\{u_2,\ldots,u_k\}\big{\}}$, so an optimal burning sequence in $G_2$ is $(u_2,u_3,\ldots,u_{k-1},u_1)$ and a minimum lazy burning set is $\{u_2,\ldots,u_{k-1}\}$. Hence, we have an infinite family of hypergraphs for which the inequality in Lemma \ref{strvswk} is strict.


We now give sufficient conditions for when a weak induced subhypergraph has (lazy) burning number no larger than that of its parent hypergraph. In fact, these conditions are the same for the lazy and round-based cases.

\begin{lemma}
\label{weak_lazy_lemma}
If $G$ is a weak induced subhypergraph of $H$ such that $|E(G)|=|E(H)|$ and $E(G)$ contains no singleton edges, then $b_L(G)\leq b_L(H)$.
\end{lemma}

\begin{proof}
Let $S$ be a minimum lazy burning set for $H$. Since every edge in $H$ contains at least two vertices that also belong to $G$, no fire can propagate in $H$ unless a vertex of $G$ is in the lazy burning set. Hence, $S\cap V(G) \neq \emptyset$. Let $S^\prime= S\cap V(G)$. We claim $S^\prime$ is a lazy burning set for $G$.

Consider the set  $V(G)\setminus S^\prime$ of vertices in $G$ that were not part of the lazy burning set. Label these vertices in chronological order with respect to the time step in which they catch fire in $LBG(H,S)$. That is, label $V(G)\setminus S^\prime$ as $u_1,\ldots,u_k$ such that if $i< j$ then $u_i$ catches fire in the same time step or in an earlier time step than $u_j$ in $LBG(H,S)$. Note that vertices of $V(G)\setminus S^\prime$ that catch fire in the same time step in $LBG(H,S)$ may be listed in any order relative to one another.

Suppose that $S^\prime$ is not a lazy burning set for $G$. Then, in $LBG(G,S^\prime)$, some vertex in $V(G)\setminus S^\prime$ does not catch fire through propagation by the end of the process. Let $u_q$ be the lowest-indexed vertex in $V(G)\setminus S^\prime$ that never catches fire. Then, all the vertices in $V(G)\setminus S^\prime$ that catch fire at a strictly earlier time step than $u_q$ in $LBG(H,S)$ will catch fire in $LBG(G,S^\prime)$. Without loss of generality, these vertices are $u_1,\ldots,u_{q-1}$.

Now, consider how $u_q$ catches fire in $LBG(H,S)$. Eventually, $u_q$ is part of some edge $e$ such that all other vertices in $e$ are on fire. But $e$ contains other vertices in $G$ apart from $u_q$, all of which are on fire (they were either in the lazy burning set or caught fire strictly earlier than $u_q$). Hence, in the corresponding edge $e\cap V(G)$ in $G$, the vertices apart from $u_q$ are either in $S^\prime$ or of the form $u_p$ with $p<q$. 

Therefore, in $LBG(G,S^\prime)$, eventually all vertices in $(e\cap V(G))\setminus \{u_q\}$ will be on fire. Thus, $u_q$ catches fire through propagation due to the edge $e\cap V(G)$. This is a contradiction, so $S^\prime$ is a lazy burning set for $G$. Thus, $b_L(G)\leq |S^\prime|\leq |S|=b_L(H)$.
\end{proof}

We must introduce some more notation for use in the following proof. In the round-based version of the game, denote the set of vertices that are on fire in a hypergraph $H$ at the end of round $r$ when burning according to the sequence $S$ by $F_r(H,S)$. 

\begin{lemma}
\label{weak_round_based_lemma}
If $G$ is a weak induced subhypergraph of $H$ such that $|E(G)|=|E(H)|$ and $E(G)$ contains no singleton edges, then $b(G)\leq b(H)$.
\end{lemma}

\begin{proof}
Let $S=(u_1,u_2,\ldots,u_{b(H)})$ be an optimal burning sequence for $H$. We construct a burning sequence $S^\prime$ for $G$ in which the $r^{th}$ source is either $u_r$ or an arbitrary vertex in $V(G)$. In particular, at each round $r\in\{1,2,\ldots,b(H)\}$, choose the $r^{th}$ source in $S^\prime$ as follows:

\begin{itemize}
\item[(1)] If $G$ became fully burned at the end of the previous round, then stop.

\item[(2)] Otherwise, if $u_r\notin V(G)$, then choose any unburned vertex in $V(G)$ as the $r^{th}$ source in $S^\prime$.

\item[(3)] Otherwise, if $u_r\in V(G)$ but $u_r$ was on fire at the end of round $r-1$, then choose any unburned vertex in $V(G)$ as the $r^{th}$ source in $S^\prime$.

\item[(4)] Otherwise, if $u_r\in V(G)$ and $u_r$ was not on fire at the end of round $r-1$, then choose $u_r$ as the $r^{th}$ source in $S^\prime$.
\end{itemize}

Clearly, by the above construction, $|S^\prime|\leq |S|$. Moreover, $|S^\prime|<|S|$ is only true if (1) occurs, in which case $S^\prime$ is indeed a burning sequence for $G$. We therefore assume that (1) does not occur, and hence $|S^\prime|=|S|$. We must show that $S^\prime$ leaves $G$ fully burned after the final round, $b(H)$. In particular, we show $F_r(H,S)\cap V(G)\subseteq F_r(G,S^\prime)$ for each round $r\in\{1,2,\ldots,b(H)\}$ by induction.

\noindent \emph{Base case.} Consider the earliest round $r_0$ at which $F_{r_0}(H,S)\cap V(G)\neq \emptyset$, since the inclusion clearly holds for all earlier rounds. Recall that every edge in $H$ contains at least two vertices of $G$. Hence, for any fire to spread to a vertex of $G$ in $BG(H,S)$, a vertex of $G$ must first be chosen as a source in $S$ at a strictly earlier round. Hence, $F_{r_0}(H,S)\cap V(G) = \{u_{r_0}\}$ (i.e.\ the first vertex of $G$ to be on fire in $BG(H,S)$ is a source). But by the above construction, in $BG(G,S^\prime)$, either $u_{r_0}$ is chosen as the $r_0^{th}$ source, or it was on fire in a strictly earlier round. Therefore, $F_{r_0}(H,S)\cap V(G)=\{u_{r_0}\}\subseteq F_{r_0}(G,S^\prime)$.

\noindent \emph{Inductive hypothesis.} Suppose $F_r(H,S)\cap V(G)\subseteq F_r(G,S^\prime)$ for some round $r\geq r_0$.

\noindent \emph{Inductive step.} By the construction of $S^\prime$, if $u_{r+1}\in V(G)$ then $u_{r+1}$ is on fire in BG$(G,S^\prime)$ at the end of round $r+1$. Consider any non-source vertex $v\in V(G)$ that catches fire through propagation in $BG(H,S)$ in round $r+1$. If $v$ was on fire in $BG(G,S^\prime)$ at an earlier round then there is nothing to show, so assume the contrary. We must show that $v$ catches fire through propagation in $BG(G,S^\prime)$ in round $r+1$. 

Consider an edge $e\in E(H)$ that caused fire to spread to $v$ in $BG(H,S)$. At the end of round $r$, each vertex in $e\setminus\{ v\}$ was on fire. Recall that $e\cap V(G)$ is a non-singleton edge in $G$. By the inductive hypothesis, $F_r(H,S)\cap V(G) \subseteq F_r(G,S^\prime)$, and in particular, $(e\cap V(G))\setminus\{v\}\subseteq F_r(H,S)\cap V(G) \subseteq F_r(G,S^\prime)$. Therefore, each vertex of $(e\cap V(G))\setminus\{v\}$ was on fire in $BG(G,S^\prime)$ at the end of round $r$. But then the edge $e\cap V(G)$ causes fire to spread to $v$ in $BG(G,S^\prime)$ in round $r+1$. 

Hence, $F_{r+1}(H,S)\cap V(G)\subseteq F_{r+1}(G,S^\prime)$. We therefore have $V(G)=F_{b(H)}(H,S)\cap V(G)\subseteq F_{b(H)}(G,S^\prime)$, so $F_{b(H)}(G,S^\prime)=V(G)$. That is, $S^\prime$ is a burning sequence for $G$, so $b(G)\leq |S^\prime| = |S|=b(H)$.
\end{proof}

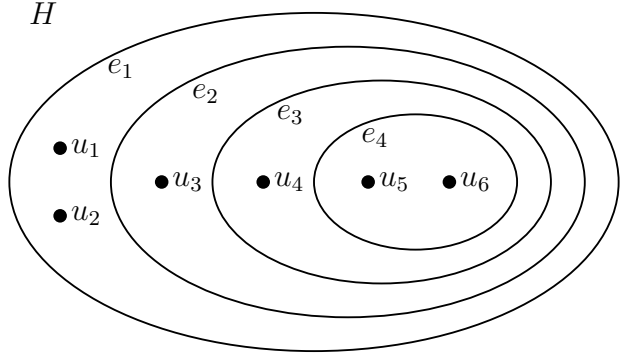
\begin{figure}
\centering
\begin{tikzpicture}[scale=0.9]
\draw [line width=0.25mm, black ] (1.5,0) ellipse (1.5cm and 1cm);
\draw [line width=0.25mm, black ] (1,0) ellipse (2.5cm and 1.5cm);
\draw [line width=0.25mm, black ] (0.5,0) ellipse (3.5cm and 2cm);
\draw [line width=0.25mm, black ] (0,0) ellipse (4.5cm and 2.5cm);
 
\node (u1) at (-4,2.5) {$H$};

\node (u1) at (-2.85,1.67) {$e_1$};
\node (u1) at (-1.6,1.31) {$e_2$};
\node (u1) at (-0.35,0.99) {$e_3$};
\node (u1) at (0.9,0.66) {$e_4$};

\node (u1) at (2,0) {};
\fill [fill=black] (u1) circle (0.1) node [right] {$u_6$};
\node (u2) at (0.8,0) {};
\fill [fill=black] (u2) circle (0.1) node [right] {$u_5$};   
\node (u3) at (-0.75,0) {};
\fill [fill=black] (u3) circle (0.1) node [right] {$u_4$};    
\node (u4) at (-2.25,0) {};
\fill [fill=black] (u4) circle (0.1) node [right] {$u_3$};    
\node (ul) at (-3.75,0.5) {};
\fill [fill=black] (ul) circle (0.1) node [right] {$u_1$}; 
\node (ul) at (-3.75,-0.5) {};
\fill [fill=black] (ul) circle (0.1) node [right] {$u_2$}; 
\end{tikzpicture}
\caption{A hypergraph $H$ that contains a weak non-induced subhypergraph with larger burning number and lazy burning number.}
\label{best_wk_non_ind}
\end{figure}

Observe that, if $H^\prime$ is a weak subhypergraph of $H$ that is not induced by any subset of $V(H)$, then $|E(H^\prime)|<|E(H)|$. The analogous results to Lemmas \ref{weak_lazy_lemma} and \ref{weak_round_based_lemma} do not hold in general if we assume the weak subhypergraphs are not induced (and hence they have fewer edges). See Figure \ref{best_wk_non_ind}, and consider the non-induced weak subhypergraph $G$ of $H$ with $V(G)=\{u_2,u_3,u_4,u_5,u_6\}$ and $E(G)=\{e_1\cap V(G)\}=\big{\{}\{u_2,u_3,u_4,u_5,u_6\} \big{\}}$. Note that $G$ is not induced by its vertex set since $E(G)\neq \{e\cap V(G)\mid e\in E(H),\ e\cap V(G)\neq\emptyset\}$ (for example, $e_2$ is not an edge in $G$). But $G$ has $b_L(G)=4>2=b_L(H)$ and $b(G)=5>4=b(H)$. To see this, first observe that $\{u_3,u_4,u_5,u_6\}$ is a minimum lazy burning set in $G$, whereas $\{u_1,u_6\}$ is a minimum lazy burning set in $H$. Also, $(u_2,u_3,u_4,u_5,u_6)$ is an optimal burning sequence in $G$, whereas $(u_6,u_4,u_1,u_2)$ is an optimal burning sequence in $H$.


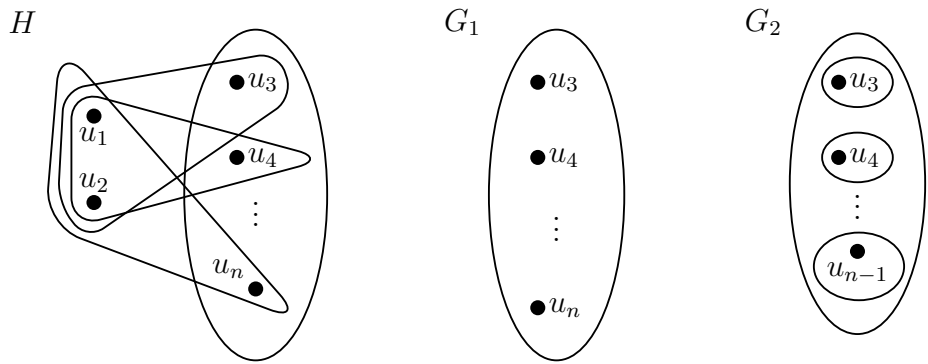
\begin{figure}
\centering
\begin{tikzpicture}[scale=1]

\node (u1) at (-6.85,2.8) {$H$};
\node (u1) at (-5.9,1.55) {};
\fill [fill=black] (u1) circle (0.1) node [below] {$u_1$};
\node (u2) at (-5.9,0.4) {};
\fill [fill=black] (u2) circle (0.1) node [above] {$u_2$};   
\node (u1) at (-4,2) {};
\fill [fill=black] (u1) circle (0.1) node [right] {$u_3$};
\node (u2) at (-4,1) {};
\fill [fill=black] (u2) circle (0.1) node [right] {$u_4$};   
\node (u3) at (-3.75,0.35) {$\vdots$};
\node (u4) at (-3.75,-0.75) {};
\fill [fill=black] (u4) circle (0.1) node [above left] {$u_n$};
\draw [line width=0.25mm, black ] (-3.75,0.5) ellipse (0.95cm and 2.2cm);    
\draw [line width=0.25mm] [black] [rounded corners=0.4cm] (-2.8,0.98)--(-6.2,0.07)--(-6.2,1.9)--cycle;
\draw [line width=0.25mm] [black] [rounded corners=0.28cm] (-3.45,2.4)--(-3.25,1.8)--(-6,-0.1)--(-6.38,0.4)--(-6.3,1.9)--cycle;
\draw [line width=0.25mm] [black] [rounded corners=0.5cm] (-3.1,-1.2)--(-6.55,0.1)--(-6.35,2.5)--cycle;

\node (u1) at (-1,2.8) {$G_1$};
\node (u1) at (0,2) {};
\fill [fill=black] (u1) circle (0.1) node [right] {$u_3$};
\node (u2) at (0,1) {};
\fill [fill=black] (u2) circle (0.1) node [right] {$u_4$};   
\node (u3) at (0.25,0.15) {$\vdots$};
\node (u4) at (0,-1) {};
\fill [fill=black] (u4) circle (0.1) node [right] {$u_n$};
\draw [line width=0.25mm, black ] (0.25,0.5) ellipse (0.9cm and 2.2cm);    

\node (u1) at (3,2.8) {$G_2$};
\node (u1) at (4,2) {};
\fill [fill=black] (u1) circle (0.1) node [right] {$u_3$};
\node (u2) at (4,1) {};
\fill [fill=black] (u2) circle (0.1) node [right] {$u_4$};   
\node (u3) at (4.25,0.45) {$\vdots$};
\node (u4) at (4.25,-0.25) {};
\fill [fill=black] (u4) circle (0.1) node [below] {$u_{n-1}$};
\draw [line width=0.25mm, black ] (4.25,0.65) ellipse (0.9cm and 2cm);    
\draw [line width=0.25mm, black ] (4.25,2) ellipse (0.47cm and 0.33cm);    
\draw [line width=0.25mm, black ] (4.25,1) ellipse (0.47cm and 0.33cm);
\draw [line width=0.25mm, black ] (4.27,-0.45) ellipse (0.6cm and 0.45cm);   

\end{tikzpicture}
\caption{$H$, $G_1$, and $G_2$ from Lemma \ref{strwk}.}
\label{HG1G2}
\end{figure}

The following result shows that weak induced subhypergraphs that do not meet the other conditions in Lemmas \ref{weak_lazy_lemma} and \ref{weak_round_based_lemma} may have larger (lazy) burning numbers than their parent hypergraphs. Indeed, in Figure \ref{HG1G2}, the weak induced subhypergraph $G_2$ of $H$ contains singleton edges, has one fewer edge than $H$, and has burning and lazy burning numbers strictly larger than those of $H$. The result also shows that strong induced subhypergraphs may have larger (lazy) burning numbers than their parent hypergraphs.

\begin{lemma}
\label{strwk}
There exist hypergraphs $H$ with strong induced subhypergraphs $G_1$ and weak induced subhypergraphs $G_2$ such that $0<b(G_1)-b(H)$, $0<b(G_2)-b(H)$, $0<b_L(G_1)-b_L(H)$, and $0<b_L(G_2)-b_L(H)$, and these differences can be arbitrarily large. 
\end{lemma}

\begin{proof}
As an example, consider $H=(V,E)$ where $V(H)=\{u_1,\ldots,u_n\}$ and $E(H)=\big{\{} \{u_1,u_2,u_3\},\{u_1,u_2,u_4\},\{u_1,u_2,u_5\},\dots,\{u_1,u_2,u_n\},\{u_3,u_4,\ldots,u_n\} \big{\}}$. Then  $(u_1,u_2,u_3)$ is an optimal burning sequence for $H$, and $\{u_1,u_2\}$ is a lazy burning set for $H$. Hence $b(H)=3$ and $b_L(H)=2$.

Now, let $G_1=(V_1,E_1)$ where $V_1=\{u_3,u_4,\ldots,u_n\}$ and $E_1=\big{\{}\{u_3,u_4,\ldots,u_n\}\big{\}}$, so $G_1$ is strongly induced by $\{u_3,u_4,\ldots,u_n\}$. Clearly, $b(G_1)=n-2$ and $b_L(G_1)=n-3$.

Finally, let $G_2=(V_2,E_2)$ where $V_2=\{u_3,u_4,\ldots,u_{n-1}\}$ and $E_2=\big{\{} \{u_3\}, \{u_4\},\ldots, \{u_{n-1}\},$ $\{u_3,u_4,\ldots,u_{n-1}\}\big{\}}$, so $G_2$ is weakly induced by $\{u_3,u_4,\ldots,u_{n-1}\}$. Clearly, $b(G_2)=n-3$ and $b_L(G_2)=n-4$.
\end{proof}


Notice that the weak induced subhypergraph $G_2$ of $H$ in Lemma \ref{strwk} need not contain fewer edges than $H$. That is, if instead $G_2$ was induced by $\{u_3,\ldots,u_n\}$, it would still have a larger (lazy) burning number than $H$. It does however contain singleton edges, and we are unsure if this must always be the case. That is, we would like to determine if there exists a class of hypergraphs $H$ with weak induced subhypergraphs $G_2$ such that $b(G_2)-b(H)$ and $b_L(G_2)-b_L(H)$ are both unbounded, and $G_2$ contains no singleton edges. Of course, in light of Lemmas \ref{weak_lazy_lemma} and \ref{weak_round_based_lemma}, such a hypergraph $G_2$ would need to have fewer edges than $H$. 

When one takes a strong subhypergraph, whether it was induced by a set of vertices or not, it is possible that the burning and lazy burning numbers both increase, and it is also possible that they both decrease. We will show this using four examples. 

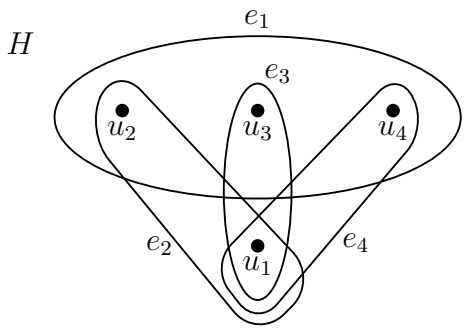
\begin{figure}
\centering
\begin{tikzpicture}[scale=0.9]

\node (u1) at (-3.5,1) {$H$};
\node (u1) at (0,1.35) {$e_1$};
\node (u1) at (-1.45,-2) {$e_2$};
\node (u1) at (0.3,0.55) {$e_3$};
\node (u1) at (1.45,-1.95) {$e_4$};

\node (u1) at (0,-2) {};
\fill [fill=black] (u1) circle (0.1) node [below] {$u_1$};
\node (u2) at (-2,0) {};
\fill [fill=black] (u2) circle (0.1) node [below] {$u_2$};   
\node (u1) at (0,0) {};
\fill [fill=black] (u1) circle (0.1) node [below] {$u_3$};
\node (u2) at (2,0) {};
\fill [fill=black] (u2) circle (0.1) node [below] {$u_4$};   

\draw [line width=0.25mm, black ] (0,-0.1) ellipse (3cm and 1.2cm);    
\draw [line width=0.25mm, black ] (0,-1.2) ellipse (0.5cm and 1.6cm);  

\draw [line width=0.25mm] [black] [rounded corners=0.5cm] (0.9,-2.5)--(0,-3.4)--(-2.55,-0.3)--(-2.1,0.7)--cycle;
\draw [line width=0.25mm] [black] [rounded corners=0.4cm] (-0.7,-2.3)--(0,-3.2)--(2.5,-0.25)--(2.1,0.6)--cycle;
\end{tikzpicture}
\caption{A hypergraph $H$ that contains an induced strong subhypergraph with larger burning number and lazy burning number.}
\label{ind_strong_ex1}
\end{figure}

For an example of an \emph{induced} strong subhypergraph $G$ with \emph{larger} burning number and lazy burning number than its parent hypergraph $H$, see Figure \ref{ind_strong_ex1}. Let $V(G)=\{u_2,u_3,u_4\}$ and $E(G)=\{e\in E(H)\mid e\subseteq V(G)\}=\{ e_1 \}$, so $G$ is a strong subhypergraph of $H$ that is induced by its vertex set. Then, $b_L(G)=2>1=b_L(H)$ and $b(G)=3>2=b(H)$. To see this, observe that $G$ is simply an edge containing three vertices. But $\{u_1\}$ is a minimum lazy burning set for $H$ and $(u_1,u_2)$ is an optimal burning sequence for $H$, so $b_L(H)=1$ and $b(H)=2$.


Also observe that this example can be extended to provide an infinite family of hypergraphs $H$ with the same property. One may construct $H$ by taking $K_{1,n}$ and adding an edge $e$ containing all $n$ of the vertices of degree one. Then, $b_L(H)=1$ and $b(H)=2$. Indeed, by taking the strong subhypergraph of $H$ induced by the vertices in $e$, we obtain a hypergraph that consists of $n$ vertices in an edge, which has lazy burning number $n-1$ and burning number $n$. Thus, the differences in (lazy) burning numbers can be arbitrarily large.

\begin{figure}[h]
\centering
\begin{tikzpicture}[scale=1.13]

\node (a) at (-3,0) {};
\fill [fill=black] (a) circle (0.105) node [right] {$u_1$};
\node (b) at (-1.5,0) {};
\fill [fill=black] (b) circle (0.105) node [left] {$u_2$};
\node (c) at (0,0) {};
\fill [fill=black] (c) circle (0.105) node [left] {$u_3$};
\node (d) at (1.5,0) {};
\fill [fill=black] (d) circle (0.105) node [left] {$u_4$};
\node (e) at (3,0) {};
\fill [fill=black] (e) circle (0.105) node [left] {$u_5$};

\draw [line width=0.25mm] [black] [rounded corners=0.85cm] (-3.7,0)--(-2.25,0.5)--(-0.8,0)--(-2.25,-0.5)--cycle;
\draw [line width=0.25mm] [black] [rounded corners=1.55cm] (-4.3,0)--(-1.6,0.9)--(1.1,0)--(-1.6,-0.9)--cycle;
\draw [line width=0.25mm] [black] [rounded corners=2.15cm] (-4.8,0)--(-1,1.3)--(2.9,0)--(-1,-1.3)--cycle;
\draw [line width=0.25mm] [black] [rounded corners=3cm] (-5.4,0)--(-0.4,1.7)--(5,0)--(-0.4,-1.7)--cycle;

\node (f) at (-2.25,0.15) {$e_1$};
\node (g) at (-0.85,0.4) {$e_2$};
\node (g) at (0.55,0.55) {$e_3$};
\node (h) at (2.2,0.65) {$e_4$};
\node (i) at (-4,1) {\textsf{$H$}};
\end{tikzpicture}

\caption{A hypergraph $H$ that contains a non-induced strong subhypergraph with larger burning number and lazy burning number.}
\label{weak_non_induced_1}
\end{figure}
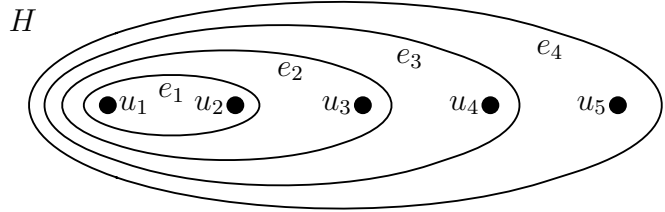

For an example of a \emph{non}-induced strong subhypergraph $G$ with \emph{larger} burning number and lazy burning number than its parent hypergraph $H$, see Figure \ref{weak_non_induced_1}. Let $V(G)=\{u_1,u_2,u_3,u_4\}$ and $E(G)=\{e_3\}$, so $G$ is a strong subhypergraph of $H$ that is not induced by its vertex set (since $e_1,e_2\notin E(G)$). Then $b_L(G)=3>1=b_L(H)$ and $b(G)=4>3=b(H)$. To see this, observe that $G$ is simply an edge containing four vertices, so $b_L(G)=3$ and $b(G)=4$. But $\{u_1\}$ is a minimum lazy burning set for $H$ and $(u_1,u_3,u_5)$ is an optimal burning sequence for $H$, so $b_L(H)=1$ and $b(H)=3$.

Again, the concept illustrated in this example can be extended to an infinite class of hypergraphs. Let the vertex set of such a hypergraph $H$ be $\{u_1,u_2,\ldots,u_n\}$, and let the edge set be $\big\{\{u_1,u_2\}, \{u_1,u_2,u_3\}, \ldots , \{u_1,\ldots, u_n\}\big\}$. Then, the hypergraph $G$ with vertex set $\{u_1,\ldots,u_{n-1}\}$ with a single edge containing all of its vertices is a non-induced strong subhypergraph of $H$. Of course, $b_L(G)=n-2$ and $b_L(H)=1$, so the differences in lazy burning numbers can be arbitrarily large. Now, observe that $b(H)=\left\lceil\frac{n+1}{2}\right\rceil$, since an optimal burning sequence for $H$ is $(u_1,u_3,\ldots,u_n)$ when $n$ is odd, and $(u_1,u_3,\ldots,u_{n-1},u_n)$ when $n$ is even. Since $b(G)=n-1$, the differences in the burning numbers of $G$ and $H$ can also be arbitrarily large.

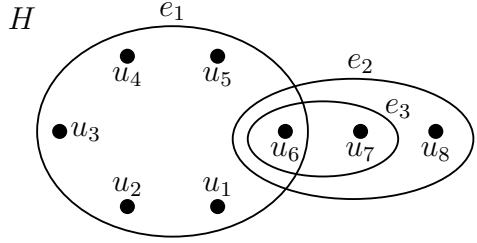
\begin{figure}
\centering
\begin{tikzpicture}[scale=1]

\node (u1) at (-3.5,1.5) {$H$};
\node (u1) at (-1.5,1.6) {$e_1$};
\node (u1) at (1,0.9) {$e_2$};
\node (u1) at (1.5,0.3) {$e_3$};

\node (u1) at (0,0) {};
\fill [fill=black] (u1) circle (0.1) node [below] {$u_6$};
\node (u2) at (-0.9,1) {};
\fill [fill=black] (u2) circle (0.1) node [below] {$u_5$};   
\node (u1) at (-2.1,1) {};
\fill [fill=black] (u1) circle (0.1) node [below] {$u_4$};
\node (u2) at (-3,0) {};
\fill [fill=black] (u2) circle (0.1) node [right] {$u_3$};   

\node (u1) at (-2.1,-1) {};
\fill [fill=black] (u1) circle (0.1) node [above] {$u_2$};
\node (u2) at (-0.9,-1) {};
\fill [fill=black] (u2) circle (0.1) node [above] {$u_1$};   
\node (u1) at (1,0) {};
\fill [fill=black] (u1) circle (0.1) node [below] {$u_7$};
\node (u2) at (2,0) {};
\fill [fill=black] (u2) circle (0.1) node [below] {$u_8$}; 

\draw [line width=0.25mm, black ] (-1.5,0) ellipse (1.8cm and 1.4cm);    
\draw [line width=0.25mm, black ] (0.5,-0.1) ellipse (1cm and 0.5cm);  
\draw [line width=0.25mm, black ] (0.9,-0.1) ellipse (1.6cm and 0.8cm); 
\end{tikzpicture}
\caption{A hypergraph $H$ that contains both an induced and a non-induced strong subhypergraph with smaller burning number and lazy burning number.}
\label{str_less_than}
\end{figure}

For an example of an \emph{induced} strong subhypergraph $G$ with \emph{smaller} burning number and lazy burning number than its parent hypergraph $H$, see Figure \ref{str_less_than}. Let $V(G)=\{u_6,u_7,u_8\}$ and $E(G)=\{e\in E(H)\mid e\subseteq V(G)\}=\{ e_2,e_3 \}$, so $G$ is a strong subhypergraph of $H$ that is induced by its vertex set. Then, $b_L(G)=1<5=b_L(H)$ and $b(G)=2<6=b(H)$. To see this, observe that $\{u_6\}$ is a minimum lazy burning set in $G$ and $(u_6,u_8)$ is an optimal burning sequence in $G$, so $b_L(G)=1$ and $b(G)=2$. But $\{u_1,u_2,u_3,u_4,u_5\}$ is a minimum lazy burning set for $H$ and $(u_6,u_1,u_2,u_3,u_4,u_5)$ is an optimal burning sequence for $H$, so $b_L(H)=5$ and $b(H)=6$.


Finally, for an example of a \emph{non}-induced strong subhypergraph $G$ with \emph{smaller} burning number and lazy burning number than its parent hypergraph $H$, see Figure \ref{str_less_than} again. This time, let $V(G)=\{u_6,u_7,u_8\}$ and $E(G)=\{ e_2 \}$, so $G$ is a strong subhypergraph of $H$ that is not induced by its vertex set (since $e_3\notin E(G)$). Then, since $G$ is simply an edge containing three vertices, and both $b_L(H)$ and $b(H)$ were found previously, it is clear that $b_L(G)=2<5=b_L(H)$ and $b(G)=3<6=b(H)$.

Both of the previous two examples (an induced/non-induced strong subhypergraph with smaller burning number) can be extended to an infinite family of hypergraphs. Indeed, one may construct $H$ with vertex set $\{u_1,\ldots,u_n\}$, and containing edges $\{u_1,\ldots,u_\ell\}$, $\{u_\ell,u_{\ell+1}\}$, \{$u_\ell,u_{\ell+1},u_{\ell+2}\}$, $\ldots$ , $\{u_\ell,\ldots,u_n\}$ for some $\ell$. Then, the induced strong subhypergraph may be obtained by simply deleting the edge $\{u_1,\ldots,u_\ell\}$, and the non-induced strong subhypergraph can be obtained by deleting all edges \emph{but} $\{u_\ell,\ldots,u_n\}$. Of course, we also delete any isolated vertices that result from this. In both cases, we can ensure the resulting subhypergraph has smaller (lazy) burning number than $H$ by making $\ell$ large enough. Informally, we can increase the (lazy) burning number of $H$ as much as we want by making sure the edge $\{u_1,\ldots,u_\ell\}$ is very large, without affecting the (lazy) burning numbers of the two subhypergraphs in question. Of course, this means that the differences in the (lazy) burning numbers between $H$ and the two subhypergraphs can be arbitrarily large.

\section{Discussion and Open Problems}

In Section \ref{general_results} we showed that no analogous bound to the Burning Number Conjecture exists in (lazy) hypergraph burning (i.e.\ a universal bound on $b_L(H)$ or $b(H)$ that is sublinear in terms of $|V(H)|$), even if we only consider hypergraphs that are both uniform and linear. We therefore want to know: what restrictions must we impose on $H$ for such a bound to exist? And what would this bound be? One possibility is that such a bound exists when we insist that $|E(H)|$ is at least $|V(H)|-1$ (which is always the case in a connected graph).

We do not know if either bound in Lemma \ref{easy_discon_lem_for_b} is tight. In Lemma \ref{sum_with_minus_k}, we improved the upper bound from Lemma \ref{easy_discon_lem_for_b}, by assuming that none of the $G_i$ are isolated vertices, and this bound is tight. It seems like the bound from Lemma \ref{sum_with_minus_k} also applies to disconnected hypergraphs with up to one isolated vertex, although a proof for this eludes us.

In Section \ref{subhypergraph_section} we proved that taking either a strong or weak induced subhypergraph can potentially increase the burning and lazy burning numbers by an arbitrarily large amount (see Lemma \ref{strwk}). However, in the only example we could find in which taking a weak induced subhypergraph increases the burning and lazy burning number, the resulting hypergraph contains singleton edges (see Figure \ref{HG1G2}). Are there examples where this does not happen? That is, does there exist a hypergraph $H$ with weak induced subhypergraph $G$ such that $G$ contains no singleton edges, $b(G)>b(H)$, and $b_L(G)>b_L(H)$? Can these differences be arbitrarily large? We also showed by examples that when one takes a strong subhypergraph, whether it was induced by a set of vertices or not, it is possible that the burning and lazy burning numbers both increase, and it is also possible that they both decrease (for a total of four different scenarios); see Figures \ref{ind_strong_ex1}, \ref{weak_non_induced_1}, and \ref{str_less_than}. We would also like to determine the necessary and sufficient conditions a hypergraph and its strong subhypergraph must satisfy in order to behave in these four ways.

We conclude by presenting additional open questions.


\begin{enumerate}
\item Given an optimal burning sequence $S=(x_1,x_2,\ldots,x_{b(H)})$, does there always exist a subset of $\{x_1,x_2,\ldots,x_{b(H)}\}$ that is a minimum lazy burning set? This is certainly true when $H$ is a graph.

\item The \emph{complement} of $H$, denoted $\overline{H}$, is the hypergraph with $V(\overline{H})=V(H)$ such that, for each subset $e$ of $V(H)$, $e$ is an edge in $\overline{H}$ if and only if $e$ is not an edge in $H$. Can the (lazy) burning number of $\overline{H}$ be bounded in terms of the (lazy) burning number of $H$? One might first investigate the case where $H$ is $k$-uniform, and only edges of size $k$ are allowed in $\overline{H}$.


\item Does the ``density'' of edges in $H$ affect $b(H)$ or $b_L(H)$? What is the best way to define density for this purpose? Intuitively, it makes sense that a hypergraph with a higher density of edges would have a lower burning number and lazy burning number. 

\item Does there exist a hypergraph (that is not a path) such that rearranging the order of an optimal burning sequence results in a non-valid burning sequence? What are the necessary and sufficient conditions for a hypergraph to have this property?

\item A hypergraph is \emph{uniquely burnable} if it has a unique optimal burning sequence up to isomorphism. What are the necessary and sufficient conditions for a hypergraph to be uniquely burnable? One might also ask the analogous question for the lazy burning game.

\item Given $n\in\mathbb{N}$, what is the set of all possible burning numbers for hypergraphs of order $n$? Is this set an interval? Can these questions be determined for certain families of hypergraphs?

\item Suppose that for any pair of edges $e_1$ and $e_2$ in a hypergraph $H$, $e_1\cap e_2\notin\{e_1,e_2\}$ (i.e. no edge in $H$ is a subset of another edge). Does this have any effect on $b(H)$ or $b_L(H)$?

\item Is there an analogous result in hypergraph burning to the Tree Reduction Theorem from graph burning? And would such a result involve hypertrees (a family of tree-like hypergraphs)?

\item Do there exist stronger bounds or even exact values of $b(H)$ and $b_L(H)$ for specific classes of hypergraphs such as projective planes, Kneser hypergraphs, and combinatorial designs? Note that the burning game on Latin square hypergraphs and Steiner triple systems has been extensively studied in \cite{Latin} and \cite{our_STS_paper} respectively.

\item Hypergraph burning using a proportion-based propagation rule is explored in \cite{our_verafin_paper}. What other propagation rules might be considered?

\item Lazy hypergraph burning is in \textbf{NP}. Is it \textbf{NP}-complete? (Note that the authors of \cite{general} have since answered this question in the affirmative).

\item In \cite{general}, the authors find several strong connections between lazy hypergraph burning and a combinatorial process known as \emph{zero-forcing}. Are there other combinatorial processes that can be linked to (lazy) hypergraph burning? 
\end{enumerate}

\section{Acknowledgements}

Authors Burgess and Pike acknowledge NSERC Discovery Grant support and Jones acknowledges support from an NSERC CGS-M scholarship and an AARMS graduate scholarship.

\bibliographystyle{abbrv}

\end{document}